\documentclass[12pt, leqno]{amsart}

\usepackage{graphicx}
\usepackage{amssymb,amsmath,amsthm,amscd}
\usepackage{mathrsfs}
\usepackage{enumerate}
\usepackage{enumitem}
\usepackage{verbatim}
\usepackage[usenames,dvipsnames]{color}
\usepackage[colorlinks=true, pdfstartview=FitV,
linkcolor=blue,citecolor=blue,urlcolor=blue]{hyperref}
\usepackage[all]{xy}
\usepackage{tikz}
\usetikzlibrary{shapes, arrows}
\usepackage{tikz-cd}
\usetikzlibrary{matrix, arrows}
\usepackage{stmaryrd} 
\usepackage{dsfont}
\usepackage{bm} 
\usepackage{scalerel}[2016/12/29] 
\usepackage{amsthm}
\usepackage{adjustbox}
\usepackage{float}
\usepackage{array}
\usepackage{multirow} 
\usepackage{caption} 
\allowdisplaybreaks[4]

\DeclareMathSymbol{\shortminus}{\mathbin}{AMSa}{"39}

\newcommand{\nc}{\newcommand}

\numberwithin{equation}{section}

\usepackage[colorinlistoftodos]{todonotes}
\usepackage{mathtools}
\mathtoolsset{showonlyrefs,showmanualtags}
\usepackage{empheq}

\theoremstyle{plain}
\newtheorem{lem}{Lemma}[section]
\newtheorem{pro}[lem]{Proposition}
\newtheorem{thm}[lem]{Theorem}
\newtheorem{cor}[lem]{Corollary}
\newtheorem{defi}[lem]{Definition}

\newcommand{\Pro}{\begin{pro}}
	\newcommand{\enpro}{\end{pro}}
\newcommand{\Lem}{\begin{lem}}
	\newcommand{\enlem}{\end{lem}}
\newcommand{\Thm}{\begin{thm}}
	\newcommand{\enthm}{\end{thm}}
\newcommand{\Cor}{\begin{cor}}
	\newcommand{\encor}{\end{cor}}
\newcommand{\Defi}{\begin{defi}}
	\newcommand{\enDefi}{\end{defi}}
\newcommand{\Proof}{\begin{proof}}
	\newcommand{\enproof}{\end{proof}}

\theoremstyle{definition} 
\newtheorem{rem}[lem]{Remark}

\newtheorem{Convention}[lem]{Convention}

\newcommand{\Conv}{\begin{Convention}}
	\newcommand{\enconv}{\end{Convention}}
\nc{\Rem}{\begin{rem}}
	\nc{\enrem}{\end{rem}}

\newcommand{\arxiv}[1]{\href{http://arxiv.org/abs/#1}{\tt arXiv:\nolinkurl{#1}}}

\nc{\epito}{\twoheadrightarrow}

\nc{\rmkend}{\hfill$\triangledown$}
\nc{\defend}{\hfill$\triangle$}

\nc{\ccc}{\mathfrak{c}}
\nc{\CCC}{\mathfrak{C}}
\nc{\Ck}{\mathfrak{C}}

\nc{\kor}{\mathbb{C}}
\nc{\indx}{\mathbb{I}}

\nc{\CC}{C}
\nc{\cc}{c}
\nc{\sss}{s}
\nc{\ck}{\mathfrak{c}}

\nc{\Bg}{B}
\nc{\Ag}{A}
\nc{\Asmone}{\pmb{\Ag}^{(-1)}(z)}
\nc{\Aps}{\pmb{\Ag}_+(z)}
\nc{\Apsone}{\pmb{\Ag}^{(1)}_+(z)}
\nc{\Apsmone}{\pmb{\Ag}^{(-1)}_+(z)}
\nc{\Ams}{\pmb{\Ag}_-(z)}
\nc{\Hg}{H}
\nc{\Thg}{\Theta}
\nc{\Thgs}{\pmb{\Theta}}
\nc{\Thgsr}{\pmb{\grave{\Theta}}}
\nc{\Apsvar}[1]{\pmb{\Ag}^{(#1)}_+(z)} 

\nc{\smin}{{\shortminus}}

\nc{\Ui}{\widetilde{\mathbf{U}}^\imath}
\nc{\DrUi}{{}^{\mathrm{Dr}}\Ui}
\nc{\fext}[2]{{#1}[\negthinspace[#2]\negthinspace]}

\nc{\gTh}{\grave{\Theta}}

\nc{\DD}{\pmb{D}}
\nc{\KK}{\mathbb{K}}

\nc{\xgp}{x^+}
\nc{\Xgp}{\pmb{x}^+}
\nc{\xgm}{x^-}
\nc{\Xgm}{\pmb{x}^-}
\nc{\xgpm}{x^{\pm}}
\nc{\psig}{\phi^+}
\nc{\phig}{\phi^-}
\nc{\phipm}{\phi^\pm}
\nc{\Psig}{\pmb{\phi^+}}
\nc{\Phig}{\pmb{\phi^-}}
\nc{\Phipm}{\pmb{\phi^\pm}}
\nc{\hg}{h}

\nc{\bY}{\mathbf{Y}}
\nc{\bT}{\mathbf{T}}

\nc{\Eg}{e}
\nc{\Kg}{K}

\nc{\Rep}{\on{Rep}}

\nc{\qq}{(q-q^{-1})^{-1}}
\nc{\factor}{\Omega}

\nc{\chring}{\Z[Y_a^{\pm 1}]_{a \in \C^\times}}
\nc{\chmap}{\chi_q}
\nc{\chmod}{\Z[\mathbf{Y}_a^{\pm 1}]_{a \in \C^\times}}
\nc{\ourchmap}{\boldsymbol{\chi}_q}
\nc{\Yring}{\mathcal{Y}} 
\nc{\Ymod}{\pmb{\mathcal{Y}}}
\nc{\Yg}{\mathbf{Y}_{i,a}}
\nc{\Ys}{\mathbf{Y}_{i,\mathbf{s}}}
\nc{\Kmat}{\mathcal{K}^0} 

\nc{\Tr}{\on{Tr}}
\nc{\id}{\on{id}}

\nc{\ourR}{\mathbf{R}}
\nc{\ourQ}{\mathbf{Q}}
\nc{\ourY}{\mathbf{Y}}

\nc{\Oq}{\mathcal{O}_q(\widehat{\g})}
\nc{\Oqi}{\mathcal{O}_q^{[i]}}
\nc{\Oqc}{\mathcal{O}_q^{\mathbf{c}}(\widehat{\g})}
\nc{\car}{\mathcal{H}}
\nc{\qaa}{U_q(\widehat{\mathfrak{g}})}
\nc{\qla}{U_q(L\mathfrak{g})}
\nc{\uqla}{\widetilde{U}_q(L\mathfrak{g})}
\nc{\drqaa}{{}^{\mathrm{Dr}}\qaa}
\nc{\uqsl}{U_q(\widehat{\mathfrak{sl}}_2)}
\nc{\uLsl}{U_qL\mathfrak{sl}_2}
\nc{\Serre}{\mathsf{Serre}}
\nc{\Sym}{\on{Sym}}
\nc{\UXp}{UX_+} 

\nc{\Tbr}{\mathbf{T}}
\nc{\degdr}{\on{deg}^{\mathrm{Dr}}}
\nc{\Uq}{\mathbf{U}}
\nc{\Uu}{\widetilde{\mathbf{U}}}
\nc{\gaf}{\widehat{\mathfrak{g}}}

\nc{\ada}{\ad_{\bhg_{-1}}}
\nc{\adb}{\ad_{\bhg_{1}}}
\nc{\adc}{\ad_{\bHg_{1}^{(2)}}}

\nc{\tAg}{\widetilde{A}}
\nc{\hAg}{\widehat{A}}

\nc{\bHg}{\overline{\Hg}}
\nc{\bhg}{\overline{\hg}}
\nc{\tHg}{\widetilde{H}}

\nc {\DrOq}{{}^{\mathrm{Dr}}\Oq}
\nc{\gup}[1]{^{(#1)}}
\nc{\gupp}{^{(1),+}}
\nc{\mi}{^{-1}}
\nc{\adh}{\operatorname{ad}_{\bar{\hg}_{-1}}}
\nc{\adhp}{\operatorname{ad}_{\bar{\hg}_{1}}}
\nc{\Omg}{\Omega^{-1}}
\nc{\Htwo}{\overline{H}_1^{(2)}}
\nc{\ad}{\operatorname{ad}}

\nc{\sspan}{\on{span}}

\newcommand{\commentout}[1]{}

\newcommand{\on}{\operatorname}
\nc{\be}{\begin{enumerate}}
	\nc{\ee}{\end{enumerate}}
\newcommand{\eq}{\begin{equation}}
	\newcommand{\eneq}{\end{equation}}
\nc{\bc}{\begin{cases}}
	\nc{\ec}{\end{cases}}
\newcommand{\eqn}{\begin{eqnarray*}}
	\newcommand{\eneqn}{\end{eqnarray*}}
\newcommand{\ba}{\begin{array}}
	\newcommand{\ea}{\end{array}}

\newcommand{\C}{{\mathbb C}}

\newcommand{\Z}{{\mathbb Z}}

\newcommand{\g}{{\mathfrak{g}}}
\nc{\Ad}{\operatorname{Ad}}

\nc{\gr}{\on{gr}}

\nc{\Aut}{\operatorname{Aut}}

\nc{\coker}{\operatorname{coker}}

\nc{\Img}{\on{Im}}
\nc{\res}{\on{res}}

\nc{\modv}[1]{{#1}\operatorname{-mod}}

\nc{\bl}{\bigl(}
\nc{\br}{\bigr)}

\newlength{\mylength}
\DeclareRobustCommand{\SkipTocEntry}[5]{}

\AtBeginDocument{%
   \def\MR#1{}
}

\title[GKLO representations for shifted QSP]
{GKLO representations for shifted quantum affine symmetric pairs}

\author[J.-R. Li]{Jian-Rong Li}
\address{Faculty of Mathematics, University of Vienna, Oskar Morgenstern Platz 1, 1090 Vienna, Austria}
\email{\href{mailto:lijr07@gmail.com}{lijr07@gmail.com}}

\author[T. Prze\'{z}dziecki]{Tomasz Prze\'{z}dziecki}
\address{Faculty of Mathematics, University of Vienna, Oskar Morgenstern Platz 1, 1090 Vienna, Austria, OrciD: 0000-0001-9700-1007}
\email{\href{mailto:tomasz.przezdziecki@univie.ac.at}{tomasz.przezdziecki@univie.ac.at}}



\thanks{The first author was supported by the Austrian Science Fund (FWF): P-34602, Grant DOI: 10.55776/P34602, and PAT 9039323, Grant-DOI 10.55776/PAT9039323, and the National Natural Science Foundation of China (No. 12471023). The second author was supported by the EPSRC grant No.\ EP/W022834/1 \emph{Kac--Moody quantum symmetric pairs, KLR algebras and generalized Schur--Weyl duality}, as well as the grant PAT 9039323, Grant-DOI 10.55776/PAT9039323.}

\begin{document}

\begin{abstract}
In this note, we introduce shifted quantum affine symmetric pairs of split simply-laced type, and construct their GKLO representations, following similar recent developments in the case of shifted twisted Yangians. A full proof that our formulas yield a representation is given. 
\end{abstract}

\maketitle

\setcounter{tocdepth}{1}
\tableofcontents


\section{Introduction}
\nc{\Uh}{U_q(\widetilde{\mathfrak{h}})}
\nc{\Uhz}{\fext{U_q(\widetilde{\mathfrak{h}})}{z}}
\nc{\ichmap}{\chi_q^\imath}

In \cite{GKLO}, Gerasimov, Kharchev, Lebedev and Oblezin (GKLO) constructed a family of remarkable infinite-dimensional representations of Yangians using explicit difference operators. These GKLO representations were later generalized to the shifted Yangian setting and used to produce quantisations of transverse slices to Schubert varieties in the affine Grassmannian \cite{KWWY,BFN}. Analogous representations for quantum affine algebras were introduced in \cite{MR2184015}, and generalized to the shifted case in \cite{FT19}, with connections to $K$-theoretic Coulomb branches and multiplicative slices. 

Recently, GKLO-type representations (called `twisted GKLO' or `$\imath$GKLO') have also been constructed for shifted twisted Yangians. More precisely, type $\mathsf{AI}$ for even weights was considered in \cite{BPT}, type $\mathsf{AIII}$ in \cite{SSX, Zichang}, and the general quasi-split case in \cite{LWW2}. Connections to the geometry of symmetric quotients of the affine Grassmannian, or fixed-point varieties under certain involutions ($\imath$slices), were investigated in \cite{BPT, LWW1}. Conjecturally, truncations of twisted Yangians defined via GKLO representations quantize the coordinate rings of such varieties. In type $\mathsf{AIII}$, it has been proven that the GKLO homomorphism factors through a Coulomb branch of cotangent type \cite{SSX, Zichang}.  More general conjectures and first results about connections to Coulomb branches of potentially non-cotangent type were also proposed in \cite{LWW1, NakajimaQSP, Sunew}. 

Even though the `rational' picture is still far from complete, in this paper we take the first step towards generalizing it to the `trigonometric' setting. Namely, we introduce shifted quantum affine symmetric pair algebras (which could also be called shifted affine $\imath$quantum groups) and construct their GKLO representations. To simplify the exposition, and be able to give a full and detailed proof, we have opted to confine ourselves to the split simply-laced type in the present paper, with more general cases to be discussed elsewhere. In analogy with the rational setting, it is natural to expect that our representations should find an interpretation in terms of $K$-theoretic Coulomb branches, a multiplicative version of $\imath$slices and $q$-W algebras of classical types. 

At the same time, the results of the present paper contribute to the ongoing programme of extending key constructions from the setting of quantum affine algebras to quantum affine symmetric pairs. The original impetus came from the discovery of Drinfeld-type presentations of quantum affine symmetric pairs \cite{bas-kol-20, lu-wang-21, LWZ-I, LWZ-II, ZhangDr, LPWZ}, followed by new developments including universal K-matrix constructions \cite{AppelVlaar, AppelVlaar4, AppelVlaar3, AppelVlaar2}, boundary $q$-characters \cite{Przez-23, LP, LP25, LP26} and $\imath$Hopf realizations \cite{iHopf}.  
The introduction of \emph{shifted} quantum affine symmetric pairs adds a new dimension to this dynamically evolving theory. 

\addtocontents{toc}{\SkipTocEntry}
	
\section*{Acknowledgements} 

We would like to thank Robin Bartlett, Lukas Tappeiner and Lewis Topley for their collaboration on a related work. 
We are also grateful to the Erwin Schr\"{o}dinger Institute (ESI) in Vienna, where parts of this project were carried out, for its hospitality and support.


\section{Shifted quantum affine symmetric pairs} 

\nc{\bA}{\mathbf{A}}
\nc{\bTh}{\boldsymbol\Theta}
\nc{\aTh}{\boldsymbol{\acute{\Theta}}}
\nc{\As}[2]{\mathbf{A}_{#1}(#2)} 
\nc{\Ths}[2]{\boldsymbol{\Theta}_{#1}(#2)} 
\nc{\aThs}[2]{\boldsymbol{\acute{\Theta}}_{#1}(#2)}
\nc{\bDel}{\boldsymbol\Delta}
\nc{\Dels}{\boldsymbol\Delta(z)} 
\nc{\ai}{\alpha^\vee_i}
\nc{\half}{\textstyle \frac12}
\nc{\bv}{\mathbf{v}}
\nc{\bu}{\mathbf{u}}
\nc{\bw}{\mathbf{w}} 
\nc{\bz}{\mathbf{z}} 
\nc{\hdiff}{D_\mu^\lambda} 
\nc{\hcor}{\C} 
\nc{\Pol}{\on{Pol}_\mu^\lambda}
\nc{\Vi}{\widetilde{\mathbb{U}}^{\imath}} 
\nc{\Ckk}{C} 
\nc{\RHS}{\on{RHS}} 
\nc{\LHS}{\on{LHS}}
\nc{\hlf}{^{\frac12}}
\nc{\mhlf}{^{-\frac12}}

Consider a finite simply-laced undirected graph $\boldsymbol\Gamma$ without edge-loops. Let $\indx$ denote its set of vertices, and let $(a_{ij})_{i,j \in \indx}$ be the corresponding generalized Cartan matrix, to which we associate Lie-theoretic data in the usual way. In particular, let $\Lambda$ (resp. $\Lambda^\vee$) denote the coweight (resp. weight) lattice, with simple roots $\alpha^\vee_i$ and fundamental weights $\omega^\vee_i$ ($i \in \indx$). We endow $\Lambda$ with the usual partial order, i.e., $\mu \leq \lambda$ if and only if $\lambda - \mu$ can be written as a sum of positive coroots with non-negative integral coefficients. Let $\Lambda_+ \subset \Lambda$ be the set of dominant coweights. 
We will use the standard notation for quantum integers and brackets, i.e., 
\[
[n] = \frac{q^n - q^{-n}}{q-q^{-1}}, \qquad [x,y]_q = xy - q yx. 
\]

For the remainder of the paper, let $\mu \in \Lambda$ be \emph{any} integral coweight and set 
\[
\mu_i = \ai(\mu). 
\]
We will use the following notation for generating series: 
\[
\As{i}{z} = \sum_{r \in \Z} A_{i,r} z^r, \qquad \Ths{i}{z} = z^{-\mu_i} +  \sum_{r > -\mu_i}(q-q\mi) \Theta_{i,r} z^r, \qquad \Dels = \sum_{r \in \Z} \Ck^r z^r. 
\]
By convention, $\Theta_{i,-\mu_i} = (q-q\mi)\mi$ and $\Theta_{i,r} = 0$ for $r < -\mu_i$. 
We will also use the normalization
\[
\aThs{i}{z} = \frac{1 - q^{-2}\Ck z^2}{1 - \Ck z^2} \Ths{i}{z}. 
\]

Below we introduce a \emph{shifted} version of a quantum symmetric pair algebra (or $\imath$quantum group) of split simply-laced type. The relations are the same as in \cite[Theorem 5.1]{lu-wang-21}, but the range of the generators is different. 

\Defi
Let $\Ui_\mu$ be the algebra generated by $A_{i,s}$ and 
$\Theta_{i,r}$, as well as central elements 
$\KK_{i}^{\pm1}$, $\Ck^{\pm1}$ $(i\in \indx, \ r > -\mu_i, \ s\in\Z)$, subject to the following relations: 
\begin{align}
\bTh_i(z) \bTh_j(w) &= \bTh_j(w) \bTh_i(z), \label{eq: drel1} \\ 
\bTh_i (z)  \bA_j(w) &= \frac{(1 -q^{-a_{ij}}zw^{-1}) (1 -q^{a_{ij}} zw \CCC)}{(1 -q^{a_{ij}}zw^{-1})(1 -q^{-a_{ij}}zw \CCC)} \bA_j(w) \bTh_i (z), \label{eq: drel2} \\ 
\bA_i(z)\bA_j(w) &= \bA_j(w)\bA_i(z), \qquad\qquad\qquad\qquad \text{ if }a_{ij}=0,  \label{eq: drel3} \\ 
(q^{-1}z -w) \bA_i(z) \bA_j(w) &= (z-q^{-1}w) \bA_j(w) \bA_i(z), \quad \quad \quad \ \text{ if }a_{ij}=-1, \label{eq: drel4}
\end{align} 
\begin{align} 
&(q^2z-w) \bA_i(z) \bA_i(w) +(q^{2}w-z) \bA_i(w) \bA_i(z) =  \\ 
&\qquad\qquad =\frac{q^{-2}}{q-q^{-1}} \KK_{i} \bDel(zw) \big( (q^2z-w)\bTh_i(w) +(q^2w-z)\bTh_i(z) \big), \label{eq: drel5}\\
&\Sym _{w_1,w_2} \big( \bA_i(w_1)\bA_i(w_2)\bA_{j}(z) -[2] \bA_i(w_1)\bA_{j}(z)\bA_i(w_2)+\bA_{j}(z)\bA_i(w_1)\bA_i(w_2) \big) = \\
&\qquad\qquad =  -\KK_{i} \frac{\bDel(w_1w_2)}{q-q^{-1}} \Sym_{w_1, w_2} \Big( [\bTh_i(w_2),\bA_j(z)]_{q^{-2}} \frac{[2] z w_1^{-1} }{1 -q^{2}w_2w_1^{-1}} \label{eq: drel6} \\
&\qquad\qquad +  [\bA_j(z),\bTh_i(w_2)]_{q^{-2}}  \frac{1 +w_2w_1^{-1}}{1 -q^{2}w_2w_1^{-1}} \Big) 
\qquad\qquad\qquad\qquad 
\text{ if }a_{ij}=-1. 
\end{align} 
\enDefi 

We will refer to \eqref{eq: drel6} as the \emph{Serre relations}. 


\section{GKLO representations} 

Fix an orientation of the graph $\boldsymbol\Gamma$. We will write $i \sim j$ if $a_{ij} = -1$ ($i$ and $j$ are connected by an edge), and $i \to j$ or $j \to i$ if there is an arrow in the given direction, depending on the orientation. 

\subsection{Representation data}

Choose a dominant integral coweight $\lambda \in \Lambda_+$ such that $\mu \leq \lambda$.   
Set 
\[ 
m_i = \omega_{i}^\vee(\lambda-\mu), \qquad 
\lambda_i = \alpha_{i}^\vee(\lambda),
\]
 $\lambda_i, m_i \in \Z_{\geq 0}$. 
In other words, 
\[
\lambda = \sum_{i \in \indx} \lambda_i  \omega_{i}, \qquad 
\mu = \sum_{i\in \indx} \mu_i \omega_{i}, \qquad 
\lambda - \mu = \sum_{i\in \indx} m_i  \alpha_{i}.  
\]  
The integers defined above satisfy the following well-known identity (see, e.g., the proof of \cite[Theorem 4.5]{KWWY}): 
\eq \label{eq: lambda mu m}
\lambda_i - \mu_i = 2m_i - \sum_{i \sim j} m_j. 
\eneq 
Additionally, we choose a tuple $\boldsymbol\theta = (\theta_i)_{i \in \indx} \in \{0,1\}^{\indx}$ satisfying the condition 
\eq \label{eq: assumption}
a_{ij} = -1 \quad \Rightarrow \quad \theta_i \theta_j = 0. 
\eneq 
We call $(\lambda, \boldsymbol\theta)$ a \emph{representation datum}. 

\Rem
Without assumption \eqref{eq: assumption}, the GKLO homomorphism defined below in Theorem \ref{thm: main GKLO} would not preserve relation \eqref{eq: drel4}. 
\enrem

\subsection{Algebra of difference operators} 
We will work in the following localized algebra of difference operators. 

\Defi 
We define the following localized algebras of polynomials and (multiplicative) difference operators. 
\begin{enumerate} 
\item 
Let $\hdiff$ be the $\C[\bz_{i,r}^{\pm \frac12}]_{i \in \indx}^{1 \leq r \leq \lambda_i}$-algebra generated by
\begin{itemize} \setlength\itemsep{0.5em}
\item $\bw_{i,k}^{\pm \frac12}, \bu_{i,k}^{\pm 1}$ (for $i \in \indx$ and $1 \leq k \leq m_i$), 
\item $(\bw_{i,k} - q^{2r-2}\bw_{i,l})^{-1}(\bw_{i,k} - q^{2r}\Ckk\mi \bw_{i,l}\mi)\mi$ (for $k \neq l$ and $r \in \Z$), 
\item $(q^{-2+2r}\bw_{i,k} - q^{-2r}C\mi \bw_{i,k}\mi)\mi$, 
\item $(q^{2r}\bw_{i,k} - C\mhlf)\mi$ if $\theta_i = 1$, 
\end{itemize}
subject to the relations
\begin{align}
\bu_{i,k}^{\epsilon_1}\bw_{j,l}^{\frac{\epsilon_2}{2}} &= q^{\epsilon_1 \epsilon_2 \delta_{ij}\delta_{kl}} \bw_{j,l}^{\frac{\epsilon_2}{2}} \bu_{i,k}^{\epsilon_1}, &\quad [\bw_{i,k}, \bw_{j,l}] &= 0 = [\bu_{i,k}, \bu_{j,l}], \\
\bu_{i,k}^{\pm1}\bu_{i,k}^{\mp1} &= 1, &\quad \bw_{i,k}^{\pm \frac12}\bw_{i,k}^{\mp \frac12} &= 1, 
\end{align}
for $\epsilon_1, \epsilon_2 \in \{\pm1\}$. 
\item Let $\Pol$ be the subalgebra of $\hdiff$ with the same generators as above \emph{except} $\bu_{i,k}^{\pm 1}$. 
\end{enumerate}
\enDefi

The algebra $\hdiff$ has a natural representation on $\Pol$  
where $\bw_{i,k}$ acts by multiplication and $\bu_{i,k}^{\pm1}$ acts by the operator sending
\[
\bw_{j,l}^n \mapsto q^{\pm n \delta_{(i,k),(j,l)}} \bw_{j,l}^n. 
\]


\subsection{Special operators} 

We will be particularly interested in operators $X_{i,k}, X'_{i,k}$ and $X''_i$ defined below. First, we need to introduce notation for certain polynomials. 
Define
\begin{align*} 
W_i^{-}(z) &= \prod_{l=1}^{m_i} (z - q^{-2}\bw_{i,l}), \qquad& 
W_i^{+}(z) &= \prod_{l=1}^{m_i} (z - \Ckk^{-1}\bw_{i,l}^{-1}), \\ 
W_{i,k}^{-}(z) &= \prod_{1 \leq k \neq l \leq m_i} (z - q^{-2}\bw_{i,l}), \qquad& 
W_{i,k}^{+}(z) &= \prod_{1 \leq k \neq l \leq m_i} (z - \Ckk^{-1}\bw_{i,l}^{-1}), \\ 
Z_i^-(z) &= \prod_{l=1}^{\lambda_i} (z - \bz_{i,l}), \qquad& Z_i^+(z)&= \prod_{l=1}^{\lambda_i} (z - C\mi \bz_{i,l}\mi)
\end{align*} 
and
\[
W_i(z) = W_i^+(z) W_i^-(z), \qquad 
W_{i,k}(z) = W_{i,k}^+(z) W_{i,k}^-(z), \qquad 
Z_i(z) =  Z_i^+(z) Z_i^-(z). 
\] 

\Defi \label{def:special-operators} 
Set 
\begin{align} 
X_{i,k} &= q \rho \bw_{i,k}\mi \prod_{l=1}^{m_i} \bw_{i,l}^{-\frac12} \left( \prod_{j \to i} \prod_{l=1}^{m_j} \bw_{j,l}^{\frac12}  \right) \left( \frac{\bw_{i,k} - qC\mhlf}{\bw_{i,k} - C\mhlf} \right)^{\theta_i} \times \\ 
&\quad \times Z_i(\bw_{i,k}) \frac{\left(\prod_{i \to j} W_j(q\mi \bw_{i,k}) \right) \left( \prod_{j \to i} W_j^+(q\mi \bw_{i,k}) \right)}
{(q^{-2}\bw_{i,k} {-} C\mi \bw_{i,k}\mi)W_{i,k}(q^{-2}\bw_{i,k})} \bu_{i,k} ,\\ 
X'_{i,k} &= - q \rho \bw_{i,k}^{-1} \prod_{l=1}^{m_i} \bw_{i,l}^{\frac12} \left( \prod_{j \to i} \prod_{l=1}^{m_j} \bw_{j,l}^{-\frac12}  \right)  \left(  \frac{q^{-2} \bw_{i,k} - q\mi C\mhlf}{q^{-2}\bw_{i,k} - C\mhlf}  \right)^{\theta_i} \times \\
&\quad \times \frac{
 \prod_{j \to i} W_j^{-}(q^{-3} \bw_{i,k})}{(q^{-2}\bw_{i,k} {-} C\mi \bw_{i,k}\mi) W_{i,k}(q^{-2}\bw_{i,k})} \bu_{i,k}\mi, \\ \label{eq: X'' defi}
 X''_i &= \frac{\prod_{l=1}^{\lambda_i} (\bz_{i,l}\hlf - C\mhlf \bz_{i,l}\mhlf)\prod_{i \sim j}\prod_{l=1}^{m_j}(\bw_{j,l}^{\frac12} - qC^{-\frac12}\bw_{j,l}^{-\frac12})}{ \prod_{l=1}^{m_i}(\bw_{i,l}^{\frac12} - C^{-\frac12}\bw_{i,l}^{-\frac12})(\bw_{i,l}^{\frac12} - q^2 C^{-\frac12}\bw_{i,l}^{-\frac12})}, 
\end{align} 
where
$\rho = \frac{1}{q-q\mi}$. 
\enDefi



\Lem \label{lem: X''}
We have 
\[
X_i'' = \eta_i \frac{Z_i^-(C\mhlf)\prod_{i \sim j}W_j^-(q\mi C\mhlf)}{W_i(C\mhlf)} 
= \eta'_i \frac{Z_i^+(C\mhlf)\prod_{i \sim j}W_j^+(q\mi C\mhlf)}{W_i(q^{-2}C\mhlf)}, 
\]
for 
\begin{align*}
\eta_i &= (-1)^{\mu_i - m_i} q^{-2m_i+2 \sum_{i \sim j}m_j} C^{\frac{-m_i}{2}} \prod_{i \sim j} \prod_{l=1}^{m_j}  \bw_{j,l}^{-\frac{1}{2}} \prod_{l=1}^{\lambda_i} \bz_{i,l}\mhlf, \\ 
\eta'_i &= (-1)^{m_i} q^{-4m_i+\sum_{i \sim j}m_j} C^{\frac{\lambda_i - m_i + \sum_{i \sim j} m_j}{2}} \prod_{i \sim j} \prod_{l=1}^{m_j}  \bw_{j,l}^{\frac{1}{2}} \prod_{l=1}^{\lambda_i} \bz_{i,l}\hlf.  
\end{align*} 
Hence
\[
(X''_i)^2 = (-C\hlf)^{\mu_i} q^{3(\mu_i - \lambda_i)} \frac{Z_i(C\mhlf)\prod_{i \sim j}W_j(q\mi C\mhlf)}{W_i(C\mhlf)W_i(q^{-2}C\mhlf)}. 
\]
\enlem 

\Proof 
First consider the numerator in \eqref{eq: X'' defi}. Since 
\[ \bz_i^{\hlf}-C^{\mhlf}\bz_i^{\mhlf} = -\bz_i^{\mhlf}(C^{\mhlf}-\bz_i), \qquad 
q^{-1}C^{\mhlf}-q^{-2}\bw_{j,r} = -q^{-2} (\bw_{j,r}-qC^{\mhlf}), 
\]
we get 
\begin{align*}
\prod_{r=1}^{\lambda_i} \left(\bz_{i,r}^{\hlf}-C^{\mhlf}\bz_{i,r}^{\mhlf}\right) &= (-1)^{\lambda_i} \left(\prod_{r=1}^{\lambda_i}\bz_{i,r}^{\mhlf}\right) Z_i^-(C^{\mhlf}), \\
 \prod_{r=1}^{m_j} \left(\bw_{j,r}^{\hlf} -qC^{\mhlf}\bw_{j,r}^{\mhlf}\right) &= (-1)^{m_j}q^{2m_j} \left(\prod_{r=1}^{m_j}\bw_{j,r}^{\mhlf}\right) W_j^-(q^{-1}C^{\mhlf}).
\end{align*}
For the denominator, since 
\[ 
\left(\bw_{i,r}^{\hlf}-C^{\mhlf}\bw_{i,r}^{\mhlf}\right) \left(\bw_{i,r}^{\hlf}-q^2C^{\mhlf}\bw_{i,r}^{\mhlf}\right) = \bw_{i,r}^{-1} (\bw_{i,r}-C^{\mhlf})(\bw_{i,r}-q^2C^{\mhlf}), 
\] 
we get 
\begin{align*} 
\prod_{r=1}^{m_i} \left(\bw_{i,r}^{\hlf} -C^{\mhlf}\bw_{i,r}^{\mhlf}\right) \left(\bw_{i,r}^{\hlf} -q^2C^{\mhlf}\bw_{i,r}^{\mhlf}\right)
= (-1)^{m_i}C^{\frac{m_i}{2}}q^{2m_i} W_i(C^{\mhlf}). 
\end{align*}
Substituting these identities into the definition of $X_i''$ and using \eqref{eq: lambda mu m} gives the first equality of the Lemma. The proof of the second identity is similar. Finally, the formula for $(X''_i)^2$ results from the multiplication of the two formulae for $X''_i$ just established. 
\enproof

\subsection{GKLO homomorphism} 

The following theorem is the main result of the paper. Its rational (i.e., twisted Yangian) version was obtained by the second author with Bartlett and Tappeiner \cite{BPT}, and independently by Lu, Wang and Weekes \cite{LWW2} (also  see \cite{LWW2, SSX, Zichang} for extensions to the quasi-split case). Our construction in \cite{BPT} was motivated by the realization in terms of Sklyanin minors, following a similar strategy in the original paper  by Gerasimov, Kharchev, Lebedev and Oblezin \cite{GKLO}. The (non-Cartan) generating series $\mathbf{b}_i(z)$ act by linear combinations of shift and inverse shift operators, with the coefficients largely constrained by the defining relations. The quantum symmetric pair version considered here can be regarded as a quantization of the GKLO formulae for twisted Yangians. The polynomials $Z_i(z)$ and $W_i(z)$ are indeed trigonometric analogues of the polynomials featuring in the latter construction. However, this ``naive'' quantization is not sufficient as the quantum symmetric pair relations force the introduction of additional factors in the GKLO formulae (e.g., the factor $\prod_{l=1}^{m_i} \bw_{i,l}^{-\frac12} \left( \prod_{j \to i} \prod_{l=1}^{m_j} \bw_{j,l}^{\frac12}  \right)$ in the formula for $X_{i,k}$). 

\Thm \label{thm: main GKLO} 
Choose a representation datum $(\lambda, \boldsymbol\theta)$ as well as nonzero complex numbers $C$ and $(\kappa_i)_{i \in \indx}$. 
There exists a unique algebra homomorphism 
\eq \label{eq: GKLO homo}
\Phi_\mu^\lambda \colon \Ui_\mu[\bz_{i,r}^{\pm \frac12}]_{i \in \indx}^{1 \leq r \leq \lambda_i} \ \longrightarrow \ \hdiff 
\eneq
given by: 
\begin{align} 
\KK_i &\mapsto \kappa_i \in \C^\times, \qquad \CCC \mapsto C \in \C^\times, \\ \label{eq: Theta action}
\ \ \ \ \aThs{i}{z} &\mapsto   \rho'_i \left( \frac{(z-qC\mhlf)(z-q\mi C\mhlf)}{(z-C\mhlf)^2}  \right)^{\theta_i}
 \frac{Z_i(z)\prod_{j \sim i} W_j(q^{-1}z)}{z^{\mu_i}W_i(z) W_i(q^{-2}z)} =: \vartheta_i(z), \\ 
\As{i}{z} &\mapsto \theta_i \tau_i \delta(z C^{\frac12}) X''_i + \sum_{k=1}^{m_i} z^{-\mu_i} \kappa_i \rho'_i \delta(z/\bw_{i,k}) X_{i,k} + \sum_{k=1}^{m_i}  \delta(z q^{-2}\Ckk \bw_{i,k}) X'_{i,k}, 
\end{align}
where by the RHS of \eqref{eq: Theta action} we mean the expansion of the rational function $\vartheta_i(z)$ in powers of $z$. The normalizing constants $\rho'_i$ and $\tau_i$ are given by
\[ 
\rho'_i = \frac{C^{\mu_i}}{q^{2(\lambda_i - \mu_i)}}, \qquad 
\tau_i = {\left( \frac{(-1)^{\mu_i}\kappa_i \rho'_i q^{3(\lambda_i - \mu_i)}}{(1+q)^2}   \right)}\hlf. 
\]
\enthm 

We will prove Theorem \ref{thm: main GKLO} in \S \ref{sec: non-Serre}--\ref{sec: Serre} below. Section \ref{sec: Serre} contains the proof that $\Phi_\mu^\lambda$ preserves the Serre relations. All the other relations are established in \S \ref{sec: non-Serre}. 

\Rem 
In the unshifted case, the generating series $\Ths{i}{z}$ arise naturally in the $\imath$Hall algebra construction of the affine $\imath$quantum group (see, e.g., \cite{lu-ruan-wang-23}). 
On the other hand, the generators $\aThs{i}{z}$ appeared earlier in \cite{bas-kol-20} in the rank one case (the $q$-Onsager algebra). It was shown in \cite{Przez-23} that, from an algebraic point of view, the series $\aThs{i}{z}$ enjoy better properties - for example, they admit a simpler coproduct formula (see \cite[Corollary 4.11]{Przez-23}). The GKLO homomorphism in Theorem \ref{thm: main GKLO} is also easier to state in terms of the series $\aThs{i}{z}$. 
\enrem

Composing the GKLO homomorphism \eqref{eq: GKLO homo} with the natural representation of $\hdiff$ on $\Pol$, we obtain the corresponding \emph{GKLO representation} of $\Ui_\mu[\bz_{i,r}^{\pm \frac12}]$. We also call the image of $\Ui_\mu[\bz_{i,r}^{\pm \frac12}]$ in $\hdiff$ under $\Phi_\mu^\lambda$ the $\lambda$\emph{-truncation} of $\Ui_\mu[\bz_{i,r}^{\pm \frac12}]$. 

\Rem
The algebra $\hdiff$ is essentially a localized quantum torus. Realizations of quantum algebras in terms of the quantum torus have been of great interest, especially from the point of view of cluster theory. In the case of quantum symmetric pairs, such a realization was recently proposed in \cite{MR4967120}. 
Theorem \ref{thm: main GKLO} gives a new alternative realization. 
\enrem 

\Rem
The formulae in Theorem \ref{thm: main GKLO} bear an obvious resemblance to Gelfand--Tsetlin formulas in the orthogonal case (see, e.g., \cite{GavKl} and \cite[Theorems 4.3, 6.1]{LP} for the quantum symmetric pair setting). We expect a more precise relationship between the two can be established via parabolic Verma modules, as in \cite[Proposition 2.175]{FPT}. 
\enrem 

\Rem
Our GKLO formulae exhibit a striking asymmetry in terms of the form of the operators $X_{i,k}$ and $X'_{i,k}$. A geometric interpretation of this asymmetry is that truncated quantum symmetric pairs or twisted Yangians of split type do not appear to have a realization as Coulomb branches of cotangent type as in \cite{BFN}. Recently, an alternative  construction via Coulomb branches of \emph{separated} cotangent type was proposed in \cite{Sunew}. 
\enrem 

\Rem
Shortly after the first version of this paper was published, a similar construction of GKLO representations, also including the quasi-split case, appeared in \cite{LWW3}. Up to normalization, the two constructions and the associated formulae appear the same. 
\enrem

\section{Proof: part I} 

In this section, we show that the GKLO homomorphism preserves all the defining relations except the Serre relation. We begin by proving some auxiliary lemmas. 

\label{sec: non-Serre}

\subsection{Auxiliary results} 

We first note that the rational function $\vartheta_i(z)$ satisfies a ``$C$-symmetry property'', analogous to that in \cite[Proposition 3.1]{Przez-23}. 

\Lem \label{lem: C-inv}
The rational function $\vartheta_i(z)$ from \eqref{eq: Theta action} is invariant under the transformation $z \mapsto C\mi z\mi$. Moreover, $\vartheta_i(z) \in z^{-\mu_i} + z^{-\mu_i+1}\C[[z]]$ as a series in $z$. 
\enlem 

\Proof
The invariance follows by a straightforward calculation. It is also clear from the formula \eqref{eq: Theta action} (thanks to the factor $z^{\mu_i}$ in the denominator) that the series $\vartheta_i(z)$ has the required range of powers. It remains to show that the lowest coefficient is $1$. Indeed, we get
\[
z^{\mu_i}\vartheta_i(z)|_{z=0} = \rho'_i \frac{Z_i(0)\prod_{j\sim i}W_j(0)}{W_i(0)^2} = \rho'_i C^{2m_i-\lambda_i - \sum_{j \sim i} m_j} q^{4m_i-2\sum_{j \sim i} m_j} = \rho'_i \frac{q^{2(\lambda_i - \mu_i)}}{C^{\mu_i}} = 1,  
\]
where the penultimate equality follows from \eqref{eq: lambda mu m}. 
\enproof

\Lem \label{lem: XjTi}
We have 
\begin{align}
X_{j,k}\aThs{i}{z} = \alpha_{i,j,k} \aThs{i}{z} X_{j,k}, \qquad 
X'_{j,k}\aThs{i}{z} = \alpha'_{i,j,k} \aThs{i}{z} X'_{j,k}, 
\end{align}
where 
\begin{align}
\alpha_{i,j,k} = \frac{(z-q^{-a_{ij}}\bw_{j,k})(z-q^{a_{ij}}\Ckk\mi \bw_{j,k}\mi)}{(z-q^{a_{ij}}\bw_{j,k})(z-q^{-a_{ij}}\Ckk\mi \bw_{j,k}\mi)}, \quad \alpha'_{i,j,k} = \frac{(z-q^{a_{ij}-2}\bw_{j,k})(z-q^{2-a_{ij}}\Ckk\mi \bw_{j,k}\mi)}{(z-q^{-2-a_{ij}}\bw_{j,k})(z-q^{2+a_{ij}}\Ckk\mi \bw_{j,k}\mi)}. 
\end{align}
Moreover, 
\[
\delta(w/\bw_{j,k}) \alpha_{i,j,k} = \delta(w/\bw_{j,k}) \alpha_{i,j}, \qquad \delta(w q^{-2}\Ckk \bw_{j,k})  \alpha'_{i,j,k} = \delta(w q^{-2}\Ckk \bw_{j,k})  \alpha_{i,j}, 
\]
where 
\[
\alpha_{i,j} = \frac{(z-q^{-a_{ij}}w)(z-q^{a_{ij}}\Ckk\mi w\mi)}{(z-q^{a_{ij}}w)(z-q^{-a_{ij}}\Ckk\mi w\mi)}. 
\]
Finally,
\[
\delta(w C^{\frac12}) X''_j \aThs{i}{z} = \delta(w C^{\frac12}) \alpha_{i,j} \aThs{i}{z} X''_j. 
\]
\enlem

\Proof
Clearly, the rational fraction coefficients on $X_{j,k}$, $X'_{j,k}$ and $X''_j$ are irrelevant, so we may just as well assume $X_{j,k} = \bu_{j,k}$, $X'_{j,k} = \bu_{j,k}\mi$ and $X''_j = 1$. Moreover, the factor $( \cdots )^{\theta_i}$ in $\aThs{i}{z}$ commutes with both $X_{j,k}$ and $X'_{j,k}$, so we may, without loss of generality, assume that $\theta_i = 0$. 
 First, let $a_{ij} = -1$. Then 
\begin{align*}
X_{j,k} \aThs{i}{z} &= \bu_{j,k} Y Z_i(z) \frac{\prod_{s \sim i} W_s(q^{-1}z)}{W_i(z) W_i(q^{-2}z)} \\ 
&= Y Z_i(z) \frac{q^{-2}(z - q\bw_{j,k})(z - q\mi C\mi \bw_{j,k}\mi) W_{j,k}(q\mi z) \prod_{j\neq s \sim i} W_s(q^{-1}z)}{W_i(z) W_i(q^{-2}z)} \bu_{j,k}, \\
 \aThs{i}{z} X_{j,k} &= Y Z_i(z) \frac{q^{-2}(z - q^{-1}\bw_{j,k})(z - qC\mi \bw_{j,k}\mi) W_{j,k}(q\mi z) \prod_{j\neq s \sim i} W_s(q^{-1}z)}{W_i(z) W_i(q^{-2}z)} \bu_{j,k},
\end{align*}
for $Y = \rho_i' z^{-\mu_i}$, 
so 
\[
X_{j,k}\aThs{i}{z} =  \frac{(z-q\bw_{j,k})(z-q^{-1}\Ckk\mi \bw_{j,k}\mi)}{(z-q^{-1}\bw_{j,k})(z-q\Ckk\mi \bw_{j,k}\mi)} \aThs{i}{z} X_{j,k}. 
\] 

Next, let $a_{ij} = 2$, i.e., $i=j$. Then 
\begin{align*}
X_{i,k} \aThs{i}{z} &=  \bu_{i,k} Y Z_i(z) \frac{\prod_{s \sim i} W_s(q^{-1}z)}{W_{i}(z) W_{i}(q^{-2}z)} \\ 
&=  \frac{q^4 Y Z_i(z) \prod_{s \sim i} W_s(q^{-1}z)}{(z{-}\bw_{i,k})(z{-}q^{-2}C\mi \bw_{i,k}\mi)(z{-}q^2\bw_{i,k})(z{-}C\mi \bw_{i,k}\mi)W_{i,k}(z) W_{i,k}(q^{-2}z)} \bu_{i,k},  \\
 \aThs{i}{z} X_{i,k} &=   \frac{q^4 Y Z_i(z) \prod_{s \sim i} W_s(q^{-1}z)}{(z{-}q^{-2}\bw_{i,k})(z{-}C\mi \bw_{i,k}\mi)(z{-}\bw_{i,k})(z{-}q^2C\mi \bw_{i,k}\mi)W_{i,k}(z) W_{i,k}(q^{-2}z)} \bu_{i,k},
\end{align*}
for $Y = \rho_i' z^{-\mu_i}$, 
so
\[
X_{i,k}\aThs{i}{z} =  \frac{(z-q^{-2}\bw_{i,k})(z-q^{2}\Ckk\mi \bw_{i,k}\mi)}{(z-q^{2}\bw_{i,k})(z-q^{-2}\Ckk\mi \bw_{i,k}\mi)} \aThs{i}{z} X_{i,k}. 
\] 
The calculation for $X'_{j,k}$ is similar. The remaining statements follow directly from the properties of the delta function. 

Finally, note that $X''_j \aThs{i}{z} = \aThs{i}{z} X''_j$, and $\delta(w C^{\frac12}) = \delta(w C^{\frac12}) \alpha_{i,j}$, 
which implies the last statement. 
\enproof

\Lem \label{lem: XiXj}
If $(i,k) \neq (j,l)$ then 
\[
X_{j,l} X_{i,k} = \beta_{i,k,j,l} X_{i,k} X_{j,l}, \qquad 
X'_{j,l} X'_{i,k} = \beta_{i,k,j,l}\mi X'_{i,k} X'_{j,l},  \qquad 
X'_{j,l} X_{i,k} = \beta'_{i,k,j,l} X_{i,k} X'_{j,l},  
\]
where 
\[
\beta_{i,k,j,l} = q^{a_{ij}} \frac{\bw_{i,k} - q^{-a_{ij}} \bw_{j,l}}{\bw_{i,k} - q^{a_{ij}} \bw_{j,l}}, \qquad 
\beta'_{i,k,j,l} = q^{a_{ij}} \frac{\bw_{i,k} - q^{2-a_{ij}}\Ckk\mi \bw_{j,l}\mi}{\bw_{i,k} - q^{2+a_{ij}}\Ckk\mi \bw_{j,l}\mi}. 
\] 
Moreover, 
\begin{align*}
\delta(w/\bw_{j,l}) \delta(z/\bw_{i,k}) \beta_{i,k,j,l} &= \delta(w/\bw_{j,l}) \delta(z/\bw_{i,k}) \beta_{i,j}, \\ 
\delta(wq^{-2}C\bw_{j,l}) \delta(zq^{-2}C\bw_{i,k}) \beta_{i,k,j,l} &= \delta(wq^{-2}C\bw_{j,l}) \delta(zq^{-2}C\bw_{i,k}) \beta_{i,j}\mi, \\ 
\delta(wq^{-2}C\bw_{j,l}) \delta(z/\bw_{i,k}) \beta'_{i,k,j,l} &= \delta(wq^{-2}C\bw_{j,l}) \delta(z/\bw_{i,k}) \beta_{i,j}, \\ 
\delta(w/\bw_{j,l}) \delta(zq^{-2}C\bw_{i,k}) \beta'_{j,l,i,k} &= \delta(w/\bw_{j,l}) \delta(zq^{-2}C\bw_{i,k}) \beta_{i,j}\mi, 
\end{align*}
where 
\[
\beta_{i,j} = q^{a_{ij}} \frac{z -  q^{-a_{ij}}w}{z -  q^{a_{ij}}w}. 
\]
Finally, for any $i,j$ with $\theta_j=1$, we have 
\[
X''_j X_{i,k} = q^{a_{ij}} \frac{\bw_{i,k} - q^{-a_{ij}}C^{-\frac12}}{\bw_{i,k} - q^{a_{ij}} C^{-\frac12}} X_{i,k} X''_j, \qquad 
X''_j X'_{i,k} =  q^{a_{ij}} \frac{q^2C\mi \bw_{i,k}\mi - q^{-a_{ij}}C^{-\frac12}}{q^2C\mi \bw_{i,k}\mi - q^{a_{ij}} C^{-\frac12}} X'_{i,k} X''_j
\]
and 
\begin{align}
\delta(wC\hlf)\delta(z/\bw_{i,k}) X''_j X_{i,k} &= \beta_{i,j} \delta(wC\hlf)\delta(z/\bw_{i,k}) X_{i,k} X''_j, \label{eq:X''Xbeta1} \\ 
\delta(wC\hlf)\delta(zq^{-2}C\bw_{i,k}) X''_j X'_{i,k} &= \beta_{i,j} \delta(wC\hlf)\delta(zq^{-2}C\bw_{i,k}) X'_{i,k} X''_j. \label{eq:X''Xbeta2}
\end{align}
\enlem  

\Proof 
The case $a_{ij} = 0$ is trivial. Below we compute the other cases. Throughout the proof, we will use various assumptions which simplify the exposition and do not cause any loss of generality. For example, when computing $X_{j,l} X_{i,k} = \beta_{i,k,j,l} X_{i,k} X_{j,l}$, we may ignore the factors in the formula for $X_{j,l}$ which obviously  commute with $\bu_{i,k}$, and vice versa. The same principle applies to the other cases considered. \\

\textbf{\emph{Case: $a_{ij}=-1$.}} 

Without loss of generality, we may assume there is an arrow $i \to j$. To simplify the presentation, we may also assume, without loss of generality, that $m_s = 0$ for $s \notin \{i,j\}$, and $m_i = m_j = 1$, as well as $\lambda_s = 0$ for all $s$. Accordingly, we will omit the second subscript on the variables $\bw$ and $\bu$. Then 
\begin{align*}
X_{j,l} X_{i,k} &= Y \bw_i^{\frac12} (q\mi \bw_j - \Ckk\mi \bw_i\mi) \bu_j 
 (q\mi \bw_i - q^{-2} \bw_j)(q\mi \bw_i - \Ckk\mi \bw_j\mi) \bu_i \\ 
&= Y \bw_i^{\frac12} (q\mi \bw_j - \Ckk\mi \bw_i\mi) 
(q\mi \bw_i - \bw_j)(q\mi \bw_i - q^{-2}\Ckk\mi \bw_j\mi) \bu_j \bu_i, \\
X_{i,k}X_{j,l}  &= Y (q\mi \bw_i - q^{-2} \bw_j)(q\mi \bw_i - \Ckk\mi \bw_j\mi) \bu_i 
\bw_i^{\frac12} (q\mi \bw_j - \Ckk\mi \bw_i\mi) \bu_j \\ 
&=  q Y \bw_i^{\frac12} (q\mi \bw_i - q^{-2} \bw_j)(q\mi \bw_i - \Ckk\mi \bw_j\mi) 
(q\mi \bw_j - q^{-2}\Ckk\mi \bw_i\mi) \bu_j \bu_i, 
\end{align*} 
for $Y = \frac{q^2\rho^2\bw_j^{-\frac32}\bw_i^{-\frac32}}{(\bw_j - \Ckk\mi \bw_j\mi)(\bw_i - \Ckk\mi \bw_i\mi)}$. Hence
\[
X_{j,l} X_{i,k} = \frac{   ( \bw_{i,k} - q \bw_{j,l} )    }{   q (\bw_{i,k} - q^{-1} \bw_{j,l})   } 
X_{i,k}X_{j,l}.  
\]

Next, 
\begin{align*}
X'_{j,l} X'_{i,k} &= Y \bw_i^{-\frac12} (q^{-3} \bw_j - q^{-2} \bw_i)  \bu_j\mi \bu_i\mi, \\ 
X'_{i,k}X'_{j,l}  &= Y \bu_i\mi  \bw_i^{-\frac12} (q^{-3} \bw_j - q^{-2} \bw_i)  \bu_j\mi  \\
&= q Y \bw_i^{-\frac12} (q^{-3} \bw_j - q^{-4} \bw_i)  \bu_j\mi \bu_i\mi, 
\end{align*} 
for $Y = \frac{q^2\rho^2 \bw_j^{-\frac12}\bw_i^{-\frac12}}{(q^{-2}\bw_j - \Ckk\mi \bw_j\mi)(q^{-2}\bw_i - \Ckk\mi \bw_i\mi)}$. Hence
\[
X'_{j,l} X'_{i,k} = \frac{ q( \bw_{i,k} - q^{-1} \bw_{j,l})  }{( \bw_{i,k} - q \bw_{j,l})  }
X'_{i,k}X'_{j,l}.  
\]

Next, 
\begin{align*}
X_{j,l} X'_{i,k} &= Y \bw_i^{\frac12} (q\mi \bw_j - \Ckk\mi \bw_i\mi)  \bu_j \bu_i\mi, \\ 
X'_{i,k}X_{j,l}  &= Y \bu_i\mi\bw_i^{\frac12} (q\mi \bw_j - \Ckk\mi \bw_i\mi)  \bu_j , \\
&= q^{-1} Y \bw_i^{\frac12} (q\mi \bw_j - q^2\Ckk\mi \bw_i\mi)  \bu_j\bu_i\mi, 
\end{align*} 
for $Y = \frac{-q^2\rho^2\bw_j^{-\frac12}\bw_i^{-\frac32}}{(\bw_j - \Ckk\mi \bw_j\mi)(q^{-2}\bw_i - \Ckk\mi \bw_i\mi)}$. Hence
\[
X_{j,l} X'_{i,k} = \frac{ q ( \bw_{j,l} - q \Ckk^{-1} \bw_{i,k}^{-1} ) }{ (\bw_{j,l} - q^{3} \Ckk^{-1} \bw_{i,k}^{-1} ) }
X'_{i,k}X_{j,l}.  
\]

Next, 
\begin{align*}
X'_{j,l} X_{i,k} &= Y \bw_i^{-\frac12} (q^{-3} \bw_j - q^{-2} \bw_i)  \bu_j\mi 
 (q\mi \bw_i - q^{-2} \bw_j)(q\mi \bw_i - \Ckk\mi \bw_j\mi) \bu_i \\ 
&= Y \bw_i^{-\frac12} (q^{-3} \bw_j - q^{-2} \bw_i)  
 (q\mi \bw_i - q^{-4} \bw_j)(q\mi \bw_i - q^2\Ckk\mi \bw_j\mi) \bu_j\mi \bu_i, \\
X_{i,k}X'_{j,l}  &= Y  (q\mi \bw_i - q^{-2} \bw_j)(q\mi \bw_i - \Ckk\mi \bw_j\mi) \bu_i \bw_i^{-\frac12} (q^{-3} \bw_j - q^{-2} \bw_i)  \bu_j\mi \\
&= q\mi Y  (q\mi \bw_i - q^{-2} \bw_j)(q\mi \bw_i - \Ckk\mi \bw_j\mi) \bw_i^{-\frac12} (q^{-3} \bw_j -  \bw_i)  \bu_j\mi \bu_i,  
\end{align*} 
for $Y = \frac{-q^2\rho^2\bw_j^{-\frac32}\bw_i^{-\frac12}}{(q^{-2}\bw_j - \Ckk\mi \bw_j\mi)(\bw_i - \Ckk\mi \bw_i\mi)}$. Hence
\[
X'_{j,l} X_{i,k} =  \frac{  ( \bw_{i,k}  -  q^3  \Ckk^{-1} \bw_{j,l}^{-1} )    }{  q ( \bw_{i,k} - q \Ckk^{-1} \bw_{j,l}^{-1} ) } 
X_{i,k}X'_{j,l}.  
\]

\textbf{\emph{Case: $a_{ij}=2$.}} 

Now suppose $i=j$ and $k \neq l$. To simplify the presentation, we may assume, without loss of generality, that $m_s = 0$ for $s \neq i$, and $m_i = 2$, as well as $\lambda_s = 0$ for all $s$. Then 
\begin{align*}
X_{i,l} X_{i,k} &= Y \bw_{i,l}^{\frac12}\bw_{i,k}^{-\frac12} \big((\bw_{i,l} - \Ckk\mi \bw_{i,l}\mi) (q^{-2}\bw_{i,l} - \Ckk\mi \bw_{i,k}\mi) (q^{-2}\bw_{i,l} - q^{-2}\bw_{i,k})\big)\mi  \bu_{i,l} \\
& \quad \times  \bw_{i,k}^{\frac12}\bw_{i,l}^{-\frac12} \big((\bw_{i,k} - \Ckk\mi \bw_{i,k}\mi) (q^{-2}\bw_{i,k} - \Ckk\mi \bw_{i,l}\mi) (q^{-2}\bw_{i,k} - q^{-2}\bw_{i,l})\big)\mi \bu_{i,k} \\ 
&= Y  q\mi \big((\bw_{i,l} - \Ckk\mi \bw_{i,l}\mi) (q^{-2}\bw_{i,l} - \Ckk\mi \bw_{i,k}\mi) (q^{-2}\bw_{i,l} - q^{-2}\bw_{i,k})\big)\mi   \\
& \quad \times   \big((\bw_{i,k} - \Ckk\mi \bw_{i,k}\mi) (q^{-2}\bw_{i,k} - q^{-2}\Ckk\mi \bw_{i,l}\mi) (q^{-2}\bw_{i,k} - \bw_{i,l})\big)\mi \bu_{i,k}\bu_{i,l}, 
\end{align*}
for $Y=  q^2 \rho^2 \bw_{i,k}^{-2} \bw_{i,l}^{-2}$. 
Interchanging $k$ and $l$, we get 
\begin{align*}
X_{i,k} X_{i,l} &= Y q\mi \big((\bw_{i,k} - \Ckk\mi \bw_{i,k}\mi) (q^{-2}\bw_{i,k} - \Ckk\mi \bw_{i,l}\mi) (q^{-2}\bw_{i,k} - q^{-2}\bw_{i,l})\big)\mi   \\
& \quad \times   \big((\bw_{i,l} - \Ckk\mi \bw_{i,l}\mi) (q^{-2}\bw_{i,l} - q^{-2}\Ckk\mi \bw_{i,k}\mi) (q^{-2}\bw_{i,l} - \bw_{i,k})\big)\mi \bu_{i,k}\bu_{i,l}. 
\end{align*} 
Hence 
\[
X_{i,l} X_{i,k} = \frac{q^2(\bw_{i,k} - q^{-2}\bw_{i,l})}{(\bw_{i,k} - q^{2}\bw_{i,l})} X_{i,k} X_{i,l}.  
\] 

Next, 
\begin{align*}
X'_{i,l} X'_{i,k} &= Y \bw_{i,k}^{\frac12} \bw_{i,l}^{-\frac12}   \big((q^{-2}\bw_{i,l} {-} q^2 \Ck\mi \bw_{i,l}\mi) (q^{-2} \bw_{i,l}{ -} q^{-2} \bw_{i,k} )  (q^{-2}\bw_{i,l} {-} C^{-1} \bw_{i,k}^{-1}) \big)\mi \bu_{i,l}\mi \\
& \quad \times  \bw_{i,k}^{-\frac12} \bw_{i,l}^{\frac12} \big({(q^{-2}\bw_{i,k} {-} q^2 \Ck\mi \bw_{i,k}\mi) (q^{-2} \bw_{i,k} {-} q^{-2} \bw_{i,l}) (q^{-2}\bw_{i,k} {-} C^{-1} \bw_{i,l}^{-1} )} \big)\mi \bu_{i,k}\mi \\
&=  Y q^{-1} \big({(q^{-2}\bw_{i,l} {-} q^2 \Ck\mi \bw_{i,l}\mi) (q^{-2} \bw_{i,l} {-} q^{-2} \bw_{i,k} )  (q^{-2}\bw_{i,l} {-} C^{-1} \bw_{i,k}^{-1})} \big)\mi \\
&\quad  \times  \big({(q^{-2}\bw_{i,k} {-} q^2 \Ck\mi \bw_{i,k}\mi) (q^{-2} \bw_{i,k} {-} q^{-4} \bw_{i,l}) (q^{-2}\bw_{i,k} {-} q^2 C^{-1} \bw_{i,l}^{-1} )} \big)\mi \bu_{i,l}\mi  \bu_{i,k}\mi,
\end{align*} 
for $Y= q^2 \rho^2$. 
Interchanging $k$ and $l$, we get 
\begin{align*}
X'_{i,k}  X'_{i,l}  &=  Y q^{-1} \big( {(q^{-2}\bw_{i,k} {-} q^2 C^{-1} \bw_{i,k}\mi) (q^{-2} \bw_{i,k} {-} q^{-2} \bw_{i,l} )  (q^{-2}\bw_{i,k} {-} C^{-1} \bw_{i,l}^{-1})} \big)\mi \\
&\quad \times  \big( {(q^{-2}\bw_{i,l} {-} q^2 C^{-1} \bw_{i,l}\mi) (q^{-2} \bw_{i,l} {-} q^{-4} \bw_{i,k}) (q^{-2}\bw_{i,l} {-} q^2 C^{-1} \bw_{i,k}^{-1} )} \big)\mi \bu_{i,k}\mi  \bu_{i,l}\mi. 
\end{align*} 
Therefore
\begin{align*}
X_{i,l}' X_{i,k}' = - \frac{  q^{-2} \bw_{i,l} - q^{-4} \bw_{i,k} }{ q^{-2} \bw_{i,k} - q^{-4} \bw_{i,l} }  X_{i,k}' X_{i,l}' =  \frac{  q^{-2} (\bw_{i,k}   -  q^2\bw_{i,l}) }{ \bw_{i,k} - q^{-2} \bw_{i,l} }  X_{i,k}' X_{i,l}'.
\end{align*}

Next, 
\begin{align*}
 X_{i,k} X_{i,l}'  &= Y \bw_{i,k}^{\frac{1}{2}} \bw_{i,l}^{{-}\frac12}   \big( {( \bw_{i,k} {-} C^{{-}1} \bw_{i,k}\mi) (q^{{-}2} \bw_{i,k} {-} q^{{-}2} \bw_{i,l} ) (q^{{-}2} \bw_{i,k} {-} C^{{-}1} \bw_{i,l}^{{-}1} ) } \big)\mi \bu_{i,k} \\
&\quad \times  \bw_{i,l}^{{-}\frac{1}{2}}  \bw_{i,k}^{\frac12}    \big( {(q^{{-}2}\bw_{i,l} {-}q^2 C^{{-}1} \bw_{i,l}\mi) (q^{{-}2} \bw_{i,l} {-} q^{{-}2} \bw_{i,k} ) (q^{{-}2}\bw_{i,l} {-} C^{{-}1} \bw_{i,k}^{{-}1} )} \big)\mi \bu_{i,l}\mi \\
&= Y q \bw_{i,l}^{{-}1} \bw_{i,k}  \big( {(\bw_{i,k} {-} C^{{-}1} \bw_{i,k}\mi) (q^{{-}2} \bw_{i,k} {-} q^{{-}2} \bw_{i,l} ) (q^{{-}2} \bw_{i,k} {-} C^{{-}1} \bw_{i,l}^{{-}1} ) } \big)\mi \\
&\quad \times  \big( {(q^{{-}2}\bw_{i,l} {-}q^2 C^{{-}1} \bw_{i,l}\mi) (q^{{-}2} \bw_{i,l} {-} \bw_{i,k} ) (q^{{-}2}\bw_{i,l} {-} q^{{-}2} C^{{-}1} \bw_{i,k}^{{-}1} )} \big)\mi  \bu_{i,k} \bu_{i,l}\mi,
\end{align*}
for $Y = -q^2 \rho^2 \bw_{i,k}^{-2}$. Moreover, 
\begin{align*}
X_{i,l}' X_{i,k} &=  Y \bw_{i,l}^{{-}\frac{1}{2}}  \bw_{i,k}^{\frac12}    \big( {(q^{{-}2}\bw_{i,l} {-}q^2 C^{{-}1} \bw_{i,l}\mi) (q^{{-}2} \bw_{i,l} {-} q^{{-}2} \bw_{i,k} ) (q^{{-}2}\bw_{i,l} {-} C^{{-}1} \bw_{i,k}^{{-}1} )} \big)\mi \bu_{i,l}\mi \\
&\quad \times  \bw_{i,k}^{\frac{1}{2}} \bw_{i,l}^{{-}\frac12}   \big( {(\bw_{i,k} {-} C^{{-}1} \bw_{i,k}\mi) (q^{{-}2} \bw_{i,k} {-} q^{{-}2} \bw_{i,l} ) (q^{{-}2} \bw_{i,k} {-} C^{{-}1} \bw_{i,l}^{{-}1} ) } \big)\mi \bu_{i,k} \\
&= Y q \bw_{i,l}^{{-}1} \bw_{i,k} \big( {(q^{{-}2}\bw_{i,l} {-}q^2 C^{{-}1} \bw_{i,l}\mi) (q^{{-}2} \bw_{i,l} {-} q^{{-}2} \bw_{i,k} ) (q^{{-}2}\bw_{i,l} {-} C^{{-}1} \bw_{i,k}^{{-}1} )} \big)\mi \\
&\quad \times  \big( {(\bw_{i,k} {-} C^{{-}1} \bw_{i,k}\mi) (q^{{-}2} \bw_{i,k} {-} q^{{-}4} \bw_{i,l} ) (q^{{-}2} \bw_{i,k} {-} q^2 C^{{-}1}  \bw_{i,l}^{{-}1} ) } \big)\mi  \bu_{i,l}\mi \bu_{i,k}. 
\end{align*}
Therefore 
\begin{align*}
X_{i,l}' X_{i,k} = q^2 \frac{ \bw_{i,k}- C^{-1}\bw_{i,l}\mi }{\bw_{i,k} - q^4 C^{-1}\bw_{i,l}\mi }X_{i,k} X_{i,l}'. 
\end{align*}
The remaining statements follow directly from the properties of the delta function. 

Finally, consider the relations involving $X''_j$ for $\theta_j=1$. Here, the coefficients on $X_{i,k}$ and $X'_{i,k}$ are irrelevant, so we may just as well assume that $X_{i,k} = \bu_{i,k}$ and $X'_{i,k} = \bu_{i,k}\mi$. Let $a_{ij} = -1$. 
We may then also assume that $m_i = 1$ and $m_s = 0$ for $s \neq i$, as well as $\lambda_s = 0$ for all $s$. Then 
\begin{align*}
X''_j X_{i,k} &= (\bw_{i,k}^{\frac12} - qC^{-\frac12}\bw_{i,k}^{-\frac12}) \bu_{i,k} = \bw_{i,k}\mhlf (\bw_{i,k} - qC^{-\frac12}) \bu_{i,k}, \\
X_{i,k} X''_j &= \bu_{i,k} (\bw_{i,k}^{\frac12} - qC^{-\frac12}\bw_{i,k}^{-\frac12}) = \bw_{i,k}\mhlf  (q\bw_{i,k} - C^{-\frac12}) \bu_{i,k}
\end{align*}
and so 
\[
X''_j X_{i,k} = q\mi \frac{\bw_{i,k} - qC^{-\frac12}}{\bw_{i,k} - q\mi C^{-\frac12}} X_{i,k} X''_j. 
\]
Similarly, 
\begin{align*}
X''_j X'_{i,k} &= (\bw_{i,k}^{\frac12} - qC^{-\frac12}\bw_{i,k}^{-\frac12}) \bu_{i,k}\mi = - q\mi C\hlf \bw_{i,k}\hlf  (q^2C\mi \bw_{i,k}\mi - qC\mhlf) \bu_{i,k}\mi, \\
X'_{i,k} X''_j &= \bu_{i,k}\mi (\bw_{i,k}^{\frac12} - qC^{-\frac12}\bw_{i,k}^{-\frac12}) = - q\mi C\hlf \bw_{i,k}\hlf  (q^3C\mi \bw_{i,k}\mi - C\mhlf) \bu_{i,k}\mi 
\end{align*}
and so 
\[
X''_j X'_{i,k} = q\mi \frac{q^2C\mi \bw_{i,k}\mi - qC^{-\frac12}}{q^2C\mi\bw_{i,k}\mi - q\mi C^{-\frac12}} X'_{i,k} X''_j. 
\]
The remaining statements follow directly from the properties of the delta function. The proof for $a_{ij} = 2$ is similar. 
\enproof

Let us abbreviate
\begin{align}
\Gamma(z) = \left(\frac{(z-qC\mhlf)(z-q\mi C\mhlf)}{(z-C\mhlf)^2}  \right)^{\theta_i}, \qquad 
\mathfrak{W}_{i,k}(z) = \Gamma(z)    \frac{ Z_i(z) \prod_{i \sim j} W_j(q\mi z) }
{W_{i,k}(z)W_{i,k}(q^{-2} z) }, 
\end{align} 
\[
\mathfrak{W}_{i}(z) = \frac{ Z_i(z) \prod_{i \sim j} W_j(q\mi z) }
{W_{i}(z)W_{i}(q^{-2} z) }, 
\]
One easily checks that 
\eq \label{eq: wikz prop}
\mathfrak{W}_{i,k}(C\mi z\mi) = (C z^2)^{-\mu_i-2} \mathfrak{W}_{i,k}(z). 
\eneq

\Lem \label{lem: XX' equal}
We have 
\begin{align*}
X_{i,k} X'_{i,k} &=  \frac{- \rho^2 q\mathfrak{W}_{i,k}(\bw_{i,k})}{\bw_{i,k}^2 (q^{-2}\bw_{i,k} {-} C\mi \bw_{i,k}\mi)(\bw_{i,k} {-} q^{-2}C\mi \bw_{i,k}\mi)} 
, \\ 
X'_{i,k} X_{i,k} &=  \frac{-\rho^2 q^5\mathfrak{W}_{i,k}(q^{-2}\bw_{i,k})}{\bw_{i,k}^2 (q^{-2}\bw_{i,k} {-} C\mi \bw_{i,k}\mi)(q^{-4}\bw_{i,k} {-} q^2C\mi \bw_{i,k}\mi)} 
. 
\end{align*} 
\enlem 

\Proof
Write
\[
X_{i,k}=A_{i,k}\bu_{i,k},
\qquad
X'_{i,k}=B_{i,k}\bu_{i,k}^{-1},
\]
where $A_{i,k}$ and $B_{i,k}$ denote the corresponding rational
coefficients in Definition~\ref{def:special-operators}.  Since
\[
\bu_{i,k}f(\bw_{i,k})
=f(q^2\bw_{i,k})\bu_{i,k},
\qquad
\bu_{i,k}^{-1}f(\bw_{i,k})
=f(q^{-2}\bw_{i,k})\bu_{i,k}^{-1},
\]
for any function $f$, we get 
\[
X_{i,k}X'_{i,k}
=A_{i,k}B_{i,k}\big|_{\bw_{i,k}\mapsto q^2\bw_{i,k}},
\qquad
X'_{i,k}X_{i,k}
=B_{i,k}A_{i,k}\big|_{\bw_{i,k}\mapsto q^{-2}\bw_{i,k}}.
\]
Let us compute $X_{i,k}X'_{i,k}$.  Under the substitution
$\bw_{i,k}\mapsto q^2\bw_{i,k}$ in $B_{i,k}$, the initial factor
contributes $-\rho\bw_{i,k}^{-1} \prod_{l=1}^{m_i}\bw_{i,l}^{\frac12}$. After multiplication by the initial factor of $A_{i,k}$, it contributes $-q\rho^2\bw_{i,k}^{-2}$. 
Moreover, the two $\theta_i$-dependent factors combine to
\begin{align*}
\left(
\frac{\bw_{i,k}-qC^{-\frac12}}{\bw_{i,k}-C^{-\frac12}}
\right)^{\theta_i}
\left(
\frac{\bw_{i,k} - q^{-1}C^{-\frac12}}{\bw_{i,k}-C^{-\frac12}}
\right)^{\theta_i} 
=\Gamma(\bw_{i,k}).
\end{align*}
For every arrow $j\to i$, the remaining factors combine as
\[
W_j^+(q^{-1}\bw_{i,k})W_j^-(q^{-1}\bw_{i,k})
=W_j(q^{-1}\bw_{i,k}),
\]
while the factors corresponding to $i\to j$ are already
$W_j(q^{-1}\bw_{i,k})$.  Thus the numerator obtained after
multiplication is
\[
\Gamma(\bw_{i,k})Z_i(\bw_{i,k})
\prod_{j\sim i}W_j(q^{-1}\bw_{i,k}).
\]
Similarly, the two $W_{i,k}$-factors in the denominator are $W_{i,k}(q^{-2}\bw_{i,k})W_{i,k}(\bw_{i,k})$, 
and the remaining two factors are
\[
(q^{-2}\bw_{i,k}-C^{-1}\bw_{i,k}^{-1})
(\bw_{i,k}-q^{-2}C^{-1}\bw_{i,k}^{-1}).
\]
Putting these contributions together and using the definition of $\mathfrak{W}_{i,k}$, we obtain the first identity. 
The proof of the second identity is analogous.  
\enproof 

Abbreviate
\begin{align}
\delta_1 &= \delta(z/\bw_{i,k}) \delta(w C \bw_{i,k}), \quad& 
\delta_2 &= \delta(q^2w/\bw_{i,k}) \delta(z q^{-2} C \bw_{i,k}), \label{eq: deltas1} \\ 
\delta_3 &= \delta(w/\bw_{i,k}) \delta(z  C \bw_{i,k}), \quad&  
\delta_4 &=  \delta(q^2z/\bw_{i,k}) \delta(w q^{-2} C \bw_{i,k}), \label{eq: deltas2}
\end{align} 
\eq
\delta_5 = \delta(zC\hlf)\delta(wC\hlf). 
\eneq
In the notation above, we suppressed the dependence on $k$, which will be clear from the context. 

\Lem \label{lem: DeltaThetaX}
We have an equality of series
\begin{align}
\rho (C\mi z\mi - z) \KK_{i} \bDel(zw) \big(\aTh_i(z) - \aTh_i(C\mi z\mi) \big) = \mathfrak{X}_i,
\end{align}
where
\begin{align*}
\mathfrak{X}_i &= \kappa_i \rho'_i \rho^2 \sum_{k=1}^{m_i} \bw_{i,k}^{-2-\mu_i}(q^{-2}\bw_{i,k} {-} C\mi \bw_{i,k}\mi)\mi   \times \\
&\quad \times \Big( -q^3(\delta_3 + \delta_1) \mathfrak{W}_{i,k}(\bw_{i,k}) + q^{7+2\mu_i}(\delta_4 + \delta_2) \mathfrak{W}_{i,k}(q^{-2}\bw_{i,k}) 
\Big) \\
& - 2 \theta_i \kappa_i \rho_i' {\textstyle \frac{1-q}{1+q}} C^{\frac{\mu_i-1}{2}} \delta_5 \mathfrak{W}_{i}(C\mhlf). 
\end{align*}
\enlem 

\Proof
For the sake of precision, let us carefully distinguish between the generating series $\aTh_i(z)$ of $\Ui_\mu$, the rational functions $\vartheta_i(z)$ in \eqref{eq: Theta action} and their series expansions $\tilde{\vartheta}_i(z)$ in $z$.  
Lemma \ref{lem: C-inv} (the invariance property) implies that $\tilde{\vartheta}_i(C\mi z\mi)$ is the expansion of $\vartheta_i(z)$ in $C\mi z\mi$. We can, therefore, use the well-known residue formula (see, e.g., \cite[Lemma C.2]{FT19}) to deduce that 
\begin{align*}
(&(Cz)\mi - z) (\tilde{\vartheta}_i(z)  -  \tilde{\vartheta}_i(C\mi z\mi))  = \rho'_i ((Cz)\mi - z) \sum_{k=1}^{m_i}  ( q^{-2} \bw_{i,k} - C^{-1} \bw_{i,k}^{-1} )\mi \times \\
&\times \bigg( \delta(q^2z/\bw_{i,k}) (q^{2}\bw_{i,k}\mi)^{\mu_i+1} 
\frac{  Z_i(q^{-2} \bw_{i,k}) \prod_{j \sim i} W_j(q^{-3} \bw_{i,k}) }{W_{i,k}(q^{-2} \bw_{i,k}) W_{i}(q^{-4} \bw_{i,k}) } \Gamma(q^{-2} \bw_{i,k}) \\ 
 &\quad +  \ \delta(z/\bw_{i,k})q^2 \bw_{i,k}^{-\mu_i-1} \frac{  Z_i(\bw_{i,k}) \prod_{j \sim i} W_j(q^{-1} \bw_{i,k}) }{ W_i(\bw_{i,k}) W_{i,k}(q^{-2} \bw_{i,k}) } \Gamma(\bw_{i,k})   \\ 
 &\quad  - \ \delta(zC\bw_{i,k})(C\bw_{i,k})^{\mu_i+1} \frac{  Z_i(C^{-1} \bw_{i,k}^{-1}) \prod_{j \sim i} W_j(q^{-1} C^{-1} \bw_{i,k}^{-1}) }{ W_{i,k}(C^{-1} \bw_{i,k}^{-1}) W_i(q^{-2} C^{-1} \bw_{i,k}^{-1}) } \Gamma(C\mi \bw_{i,k}\mi)   \\ 
 &\quad -  \    \delta(q^{-2}zC\bw_{i,k})q^2 (q^{-2}C\bw_{i,k})^{\mu_i+1} 
 \frac{Z_i(q^2 C^{-1} \bw_{i,k}^{-1}) \prod_{j \sim i} W_j( q C^{-1} \bw_{i,k}^{-1} ) }{W_i(q^2 C^{-1} \bw_{i,k}^{-1}) W_{i,k}(C^{-1} \bw_{i,k}^{-1}) } \Gamma(q^{2}C\mi \bw_{i,k}\mi)  \bigg) \\ 
 &+ 2\theta_i\rho_i' q^{-1}(1-q)^2 C^{\frac{\mu_i-1}{2}} \delta(z C\hlf)  \frac{Z_i(C\mhlf)\prod_{j \sim i} W_j(q^{-1}C\mhlf)}{W_i(C\mhlf) W_i(q^{-2}C\mhlf)}. 
\end{align*} 
Hence 
\begin{align*}
(&(Cz)\mi - z) (\aTh_i(z) - \aTh_i(C\mi z\mi))  = q\rho \rho'_i ((Cz)\mi - z) \sum_{k=1}^{m_i}  ( q^{-2} \bw_{i,k} - C^{-1} \bw_{i,k}^{-1} )\mi \times \\
&\times \bigg( -\delta(q^2z/\bw_{i,k}) (q^{2}\bw_{i,k}\mi)^{\mu_i+2} (q^{-4}\bw_{i,k} - C\mi \bw_{i,k}\mi )\mi \mathfrak{W}_{i,k}(q^{-2} \bw_{i,k})  \\ 
&\quad +  \ \delta(z/\bw_{i,k})q^2 \bw_{i,k}^{-\mu_i-2} (\bw_{i,k} - C\mi \bw_{i,k}\mi)\mi \mathfrak{W}_{i,k}(\bw_{i,k}) \\ 
&\quad +  \ \delta(zC\bw_{i,k})(C\bw_{i,k})^{\mu_i+2} (q^{-2}C\mi \bw_{i,k}\mi - q^{-2} \bw_{i,k})\mi \mathfrak{W}_{i,k}(C\mi \bw_{i,k}\mi)   \\ 
&\quad - \ \delta(q^{-2}zC\bw_{i,k})q^2 (q^{-2}C\bw_{i,k})^{\mu_i+2} (q^{2}C\mi \bw_{i,k}\mi - q^{-2} \bw_{i,k})\mi  \mathfrak{W}_{i,k}(q^2 C\mi \bw_{i,k}\mi)   \bigg) \\
&+ 2\theta_i\rho_i' q^{-1}(1-q)^2 C^{\frac{\mu_i-1}{2}} \delta(z C\hlf)  \mathfrak{W}_{i}(C\mhlf). 
\end{align*} 
Using \eqref{eq: wikz prop}, one gets 
\begin{align*}
((Cz)\mi - z) &(\aTh_i(z) - \aTh_i(C\mi z\mi))  = q^3 \rho \rho'_i \sum_{k=1}^{m_i} \bw_{i,k}^{-\mu_i-2}  ( q^{-2} \bw_{i,k} - C^{-1} \bw_{i,k}^{-1} )\mi \times \\
&\times ((Cz)\mi - z)  \bigg( \frac{\big(\delta(z/\bw_{i,k}) - \delta(zC \bw_{i,k})\big)\mathfrak{W}_{i,k}(\bw_{i,k}) }{ (\bw_{i,k} - C\mi \bw_{i,k}\mi)}\\
&\quad - \frac{\big(\delta(q^2z/\bw_{i,k}) - \delta(q^{-2}zC \bw_{i,k})\big)  q^{2\mu_i+4}\mathfrak{W}_{i,k}(q^{-2}\bw_{i,k})}{  (q^{-2}\bw_{i,k} - q^2C\mi \bw_{i,k}\mi )} \bigg) \\ 
&+ 2\theta_i\rho_i' q^{-1}(1-q)^2 C^{\frac{\mu_i-1}{2}} \delta(z C\hlf)  \mathfrak{W}_{i}(C\mhlf). 
\end{align*} 
The above simplifies to 
\begin{align*}
((Cz)\mi - z) &(\aTh_i(z) - \aTh_i(C\mi z\mi))  = q^3 \rho \rho'_i \sum_{k=1}^{m_i} \bw_{i,k}^{-\mu_i-2}  ( q^{-2} \bw_{i,k} - C^{-1} \bw_{i,k}^{-1} )\mi \times \\
&\times \bigg( -\big(\delta(z/\bw_{i,k}) + \delta(zC \bw_{i,k})\big)\mathfrak{W}_{i,k}(\bw_{i,k})\\
&\quad + \big(\delta(q^2z/\bw_{i,k}) + \delta(q^{-2}zC \bw_{i,k})\big)  q^{2\mu_i+4}\mathfrak{W}_{i,k}(q^{-2}\bw_{i,k}) \bigg) \\ 
&+ 2\theta_i\rho_i' q^{-1}(1-q)^2 C^{\frac{\mu_i-1}{2}} \delta(z C\hlf)  \mathfrak{W}_{i}(C\mhlf). 
\end{align*} 
Multiplying by $\rho \KK_{i} \bDel(zw)$ yields the result. 
\enproof

\subsection{Proof of the main theorem} 

First, note that Lemma \ref{lem: C-inv} implies that the operator by which $\aThs{i}{z}$ is acting has the correct range of powers of $z$. Moreover, it follows immediately from the definitions that relations \eqref{eq: drel1}, \eqref{eq: drel3}, as well as \eqref{eq: drel2} with $a_{ij} = 0$, are preserved. Let us show the remaining relations. 

\subsubsection{Relation \eqref{eq: drel2}} 

Abbreviate $\zeta = w^{-\mu_j} \kappa_j \rho'_j$. 
We have 
\begin{align}
\As{j}{w} \Ths{i}{z} &= \theta_j \tau_j \delta(w C^{\frac12}) X''_j  \Ths{i}{z} + \sum_{k=1}^{m_j} \Big(\zeta   \delta(w/\bw_{j,k}) X_{j,k}  + \delta(w q^{-2}\Ckk \bw_{j,k}) X'_{j,k} \Big) \Ths{i}{z}\\ 
&=  \alpha_{i,j} \theta_j \tau_j \delta(w C^{\frac12}) \Ths{i}{z} X''_j \\
 & \quad + \alpha_{i,j} \Ths{i}{z} \sum_{k=1}^{m_j} \Big( \zeta \delta(w/\bw_{j,k}) X_{j,k}  +  \delta(w q^{-2}\Ckk \bw_{j,k}) X'_{j,k} \Big) \\ 
&= \alpha_{i,j} \Ths{i}{z}\As{j}{w},
\end{align} 
which is equivalent to \eqref{eq: drel2}. The second equality above follows from Lemma~\ref{lem: XjTi}. 

\subsubsection{Relation \eqref{eq: drel4}} \label{subsub: d4}
We have $a_{ij} = -1$. By assumption \eqref{eq: assumption}, either $\theta_i$ or $\theta_j$ equals zero. Without loss of generality, we may assume that $\theta_i=0$. 
Abbreviate $\zeta = w^{-\mu_j} \kappa_j \rho'_j$ and $\zeta' = z^{-\mu_i} \kappa_i \rho'_i$. 
Then 
\begin{align*}
\As{j}{w} \As{i}{z} &= \sum_{l=1}^{m_j} \sum_{k=1}^{m_i} \Big( \zeta \zeta' \delta(\textstyle\frac{w}{\bw_{j,l}}) \delta(\textstyle\frac{z}{\bw_{i,k}}) X_{j,l} X_{i,k}
+  \delta(wq^{-2}C\bw_{j,l}) \delta(zq^{-2}C\bw_{i,k}) X'_{j,l} X'_{i,k} \\ 
& \quad + \zeta'\delta(wq^{-2}C\bw_{j,l}) \delta(z/\bw_{i,k}) X'_{j,l} X_{i,k}
+ \zeta\delta(w/\bw_{j,l}) \delta(zq^{-2}C\bw_{i,k})  X_{j,l} X'_{i,k} \Big) \\ 
& \quad + \theta_j \tau_j \delta(w C^{\frac12}) X''_j \sum_{k=1}^{m_i} \Big( \zeta' \delta(z/\bw_{i,k})  X_{i,k} 
+  \delta(zq^{-2}C\bw_{i,k})  X'_{i,k}  \Big) \\ 
&= q\mi \frac{z-qw}{z-q^{-1}w} \As{i}{z} \As{j}{w}, 
\end{align*}
by Lemma \ref{lem: XiXj}, which is equivalent to \eqref{eq: drel4}. 

\subsubsection{Relation \eqref{eq: drel5}}

Using Lemma \ref{lem: XiXj}, and arguing in the same way as in \S \ref{subsub: d4}, we deduce that 
\begin{align} \label{eq: AA XX'}
\LHS &= (q^2z-w) \bA_i(z) \bA_i(w) +(q^{2}w-z) \bA_i(w) \bA_i(z) = \\ 
&= \kappa_i\rho_i'(q^2z - w) \sum_{k=1}^{m_i} \left(z^{-\mu_i}\delta_1 X_{i,k}X'_{i,k} + w^{-\mu_i}\delta_2  X'_{i,k}X_{i,k} \right) \\ 
&+  \kappa_i\rho_i'(q^2w - z) \sum_{k=1}^{m_i} \left( w^{-\mu_i}\delta_3 X_{i,k}X'_{i,k} + z^{-\mu_i} \delta_4  X'_{i,k}X_{i,k} \right) \\ 
&+ 2 \theta_i \tau_i^2 C\mhlf (q^2-1) \delta_5 (X''_i)^2. 
\end{align}
Moreover,  Lemmas \ref{lem: X''} and \ref{lem: XX' equal} imply that 
\begin{align}
\delta_{1} (q^2z - w) z^{-\mu_i} X_{i,k}X'_{i,k} &= 
  \frac{-q^3 \rho^2 \delta_1}{\bw_{i,k}^{2+\mu_i} (q^{-2}\bw_{i,k} - C\mi \bw_{i,k}\mi)} 
\mathfrak{W}_{i,k}(\bw_{i,k}), \\ 
\delta_2 (q^2z - w) w^{-\mu_i} X'_{i,k}X_{i,k} &= 
 \frac{q^{7+2\mu_i} \rho^2 \delta_2}{\bw_{i,k}^{2+\mu_i} (q^{-2}\bw_{i,k} - C\mi \bw_{i,k}\mi)} 
\mathfrak{W}_{i,k}(q^{-2}\bw_{i,k}), \\
\delta_3 (q^2w - z) w^{-\mu_i} X_{i,k}X'_{i,k} &= \frac{-q^3 \rho^2 \delta_3}{\bw_{i,k}^{2+\mu_i} (q^{-2}\bw_{i,k} - C\mi \bw_{i,k}\mi)} 
\mathfrak{W}_{i,k}(\bw_{i,k}), \\ 
\delta_4 (q^2w - z) z^{-\mu_i} X'_{i,k}X_{i,k} &=   \frac{q^{7+2\mu_i} \rho^2 \delta_4}{\bw_{i,k}^{2+\mu_i} (q^{-2}\bw_{i,k} - C\mi \bw_{i,k}\mi)}  \mathfrak{W}_{i,k}(q^{-2}\bw_{i,k}), \\ 
\delta_5 (X''_i)^2 &= \delta_5 (-C\hlf)^{\mu_i} q^{3(\mu_i - \lambda_i)} \mathfrak{W}_i(C\mhlf). 
\end{align} 
Hence $\LHS = \mathfrak{X}_i$.

On the other hand, observe that the RHS of \eqref{eq: drel5} can be expressed as
\begin{align}
\RHS &= q^{-2}\rho \KK_{i} \bDel(zw) \big( (q^2z-w)\bTh_i(w) +(q^2w-z)\bTh_i(z) \big) \\ 
&= q^{-2}\rho \KK_{i} \bDel(zw) \big( (q^2C\mi w\mi-w)\bTh_i(w) +(q^2C\mi z\mi-z)\bTh_i(z) \big) \\ 
&= \rho \KK_{i} \bDel(zw) \big( (C\mi w\mi-w)\aTh_i(w) +(C\mi z\mi-z)\aTh_i(z) \big) \\ 
&= \rho (C\mi z\mi-z) \KK_{i} \bDel(zw) \big(\aTh_i(z) - \aTh_i(C\mi z\mi) \big).  
\end{align}
Lemma \ref{lem: DeltaThetaX} now implies that $\RHS = \mathfrak{X}_i$. Hence $\LHS = \RHS$.

\section{Proof: part II - Serre relations}

\label{sec: Serre} 

We will now complete the proof, by showing that the Serre relation holds as well. We begin by analyzing the LHS of the relation, identifying a part of the formula which must vanish under the GKLO homomorphism. We then calculate the remainder and compare it with the RHS. 

\nc{\TV}{\on{TV}}

\subsection{The LHS} 

Pick $i,j$ with $a_{ij} = -1$. 
Let us write the LHS of \eqref{eq: drel6} in the basis of ordered monomials in the difference operators $\bu_{i,k}^a,  \bu_{i,l}^b, \bu_{j,m}^c$ (for $a,b,c \in \{0, \pm1\}$) with some coefficients in the variables $\bw_{\star,\star}$. These can match monomials on the RHS of \eqref{eq: drel6} only if: whenever $a \neq 0$ or $b \neq 0$, then the corresponding coefficient vanishes. Let us abbreviate the sum on the LHS which needs to vanish as $\TV$ (to vanish). 

\subsubsection{The vanishing coefficients: $\theta_i = \theta_j = 0$} 
For $1 \leq k \neq l \leq m_i$ and $1 \leq m \leq m_j$, define 
\begin{align}
M_{k,m} &= \Sym_{w_1, w_2} [\delta(w_1/ \bw_{i,k}) X_{i,k}, [\delta(w_2/ \bw_{i,k}) X_{i,k}, \delta(z/\bw_{j,m})  X_{j,m}]_{q\mi}]_q, \\  
N_{k,l,m} &= \Sym_{w_1, w_2} [\delta(w_1/ \bw_{i,k}) X_{i,k}, [\delta(w_2/ \bw_{i,l}) X_{i,l}, \delta(z/\bw_{j,m})  X_{j,m}]_{q\mi}]_q. 
\end{align}  
Then \cite[\S C.(vii)]{FT19} implies that 
\eq \label{eq: MN vanish} 
M_{k,m} = N_{k,l,m} + N_{l,k,m} = 0. 
\eneq 
Let us generalize this fact. 

\Lem \label{lem: XXX commutators}
Let $1 \leq k \neq l \leq m_i$ and $1 \leq m \leq m_j$.  
The following identities hold: 
\begin{align}
[X_{i,k},[X_{i,k}, X'_{j,m}]_{q\mi}]_q &= 0,  \label{eq: Serre aux1} \\[1pt]
[X_{i,k},[X_{i,l}, X'_{j,m}]_{q\mi}]_q &= \gamma_{k,l,m} X_{i,k} X_{i,l}  X_{j,m}', \label{eq: Serre aux2}  \\[1pt]
[X_{i,k}, [X'_{i,l}, X_{j,m}]_{q\mi}]_q &= \gamma'_{k,l,m} X_{i,k} X'_{i,l}  X_{j,m}, \label{eq: Serre aux3}  \\[1pt]
[X'_{i,k}, [X_{i,l}, X_{j,m}]_{q\mi}]_q &= \gamma''_{k,l,m} X'_{i,k} X_{i,l}  X_{j,m}, \label{eq: Serre aux4}
\end{align}
where 
\begin{align}
\gamma_{k,l,m} &= \frac{- {\left(q - q^{-1} \right)}^2 \bw_{i,k} \bw_{i,l} \left(C^{-1}(q^3+q)  \bw_{j,m}^{-1} - \bw_{i,k} - \bw_{i,l}\right)}{ \left(\bw_{i,l} - q C^{-1} \bw_{j,m}^{-1}\right) \left(q^2 \bw_{i,l} - \bw_{i,k}\right) \left(\bw_{i,k} - q C^{-1} \bw_{j,m}^{-1}\right)}, \\ 
\gamma'_{k,l,m} &= \frac{q^3 {\left(q- q\mi \right)}^2 C^{-1} \bw_{i,k}  \bw_{i,l}^{-1} \left(q^3C^{-1} \bw_{i,l}^{-1} + q \bw_{i,k} - (q^2+1) \bw_{j,m}\right)}{\left(q^4 C^{-1} \bw_{i,l}^{-1} - \bw_{i,k}\right)\left(q^3 C^{-1} \bw_{i,l}^{-1} - \bw_{j,m}\right)  \left(q \bw_{i,k} - \bw_{j,m}\right)},  \\
\gamma''_{k,l,m} &= \frac{ q {\left(q - q\mi \right)}^2 C\mi  \bw_{i,l} \bw_{i,k}\mi \left(q^3C^{-1} \bw_{i,k}^{-1} + q\bw_{i,l}  - (q^2+1) \bw_{j,m}\right) }{\left(\bw_{i,l} - C\mi \bw_{i,k}\mi \right) \left( q^3 C\mi \bw_{i,k}\mi - \bw_{j,m} \right)\left(q \bw_{i,l} - \bw_{j,m}\right)}
\end{align}
\enlem 

\Proof
The lemma follows by a direct calculation using the definition of the operators $X_{i,k}$, $X'_{j,m}$, etc. We will only calculate \eqref{eq: Serre aux2} explicitly. 
By Lemma~\ref{lem: XiXj},
\begin{align}
X'_{j,m}X_{i,r}
 &= r_r X_{i,r}X'_{j,m},
&
 r_r&=q^{-1}\frac{\bw_{i,r}-q^3C^{-1}\bw_{j,m}^{-1}}
 {\bw_{i,r}-qC^{-1}\bw_{j,m}^{-1}},
\label{eq:note-r-r}\\
X_{i,l}X_{i,k}
 &=s_{k,l}X_{i,k}X_{i,l},
&
 s_{k,l}&=q^2
 \frac{\bw_{i,k}-q^{-2}\bw_{i,l}}
 {\bw_{i,k}-q^2\bw_{i,l}}.
\label{eq:note-s-kl}
\end{align}
Notice that $r_l$ commutes with $\bu_{i,k}$, $r_k$ commutes with 
$\bu_{i,l}$, and $s_{k,l}$ commutes with $\bu_{j,m}^{-1}$. Thus no additional
$q$-shift of these coefficients is introduced in the following reordering: 
\begin{align*}
[X_{i,k},[X_{i,l},X'_{j,m}]_{q^{-1}}]_q =\bigl(1-q^{-1} r_l-q s_{k,l}r_k+s_{k,l}r_kr_l\bigr)
 X_{i,k}X_{i,l}X'_{j,m}.
\end{align*}
Moreover, the scalar coefficient simplifies to  
\begin{align*}
1-q^{-1} r_l-q s_{k,l}r_k+s_{k,l}r_kr_l =\gamma_{k,l,m}, 
\end{align*}
as required. 
\enproof

Observe that, for $k \neq l$, one has 
\begin{align*} 
P_{k,m} :=& \ \Sym_{w_1, w_2} [\delta(w_1/ \bw_{i,k}) X_{i,k},[\delta(w_2/ \bw_{i,k}) X_{i,k}, \delta(zq^{-2}C\bw_{j,m})  X'_{j,m}]_{q\mi}]_q  \\ 
=& \  [X_{i,k},[X_{i,k}, X'_{j,m}]_{q\mi}]_q\Sym_{w_1, w_2} \delta(w_1/\bw_{i,k}) \delta(w_2/q^2\bw_{i,k})  \delta(zq^{-2}C\bw_{j,m})  , \\ 
R_{k,l,m} :=& \ \Sym_{w_1, w_2} [\delta(w_1/ \bw_{i,k}) X_{i,k},[\delta(w_2/ \bw_{i,l}) X_{i,l}, \delta(zq^{-2}C\bw_{j,m})  X'_{j,m}]_{q\mi}]_q \\ 
=& \ [X_{i,k},[X_{i,l}, X'_{j,m}]_{q\mi}]_q\Sym_{w_1, w_2}  \delta(w_1/\bw_{i,k}) \delta(w_2/\bw_{i,l})  \delta(zq^{-2}C\bw_{j,m}) , \\ 
S_{k,l,m} :=& \ \Sym_{w_1, w_2} [\delta(w_1/ \bw_{i,k}) X_{i,k},[\delta(w_2 q^{-2} C \bw_{i,l}) X'_{i,l}, \delta(z/\bw_{j,m})  X_{j,m}]_{q\mi}]_q \\ 
=& \  [X_{i,k},[X'_{i,l}, X_{j,m}]_{q\mi}]_q\Sym_{w_1, w_2}  \delta(w_1/\bw_{i,k}) \delta(w_2q^{-2}C\bw_{i,l})  \delta(z/\bw_{j,m}), \\ 
T_{k,l,m} :=& \ \Sym_{w_1, w_2} [\delta(w_1q^{-2}C\bw_{i,k}) X'_{i,k},[\delta(w_2/ \bw_{i,l}) X_{i,l}, \delta(z/\bw_{j,m})  X_{j,m}]_{q\mi}]_q \\ 
=& \ [X'_{i,k},[X_{i,l}, X_{j,m}]_{q\mi}]_q\Sym_{w_1, w_2}  \delta(w_1q^{-2}C\bw_{i,k}) \delta(w_2/\bw_{i,l})  \delta(z/\bw_{j,m}) .  
\end{align*} 

It follows directly from \eqref{eq: Serre aux1} that $P_{k,m} = 0$. Moreover, since, by \eqref{eq: Serre aux2}, 
\[
\gamma_{l,k,m} = \frac{q^2\bw_{i,l} - \bw_{i,k}}{q^2\bw_{i,k} - \bw_{i,l}} \gamma_{k,l,m}, \qquad 
X_{i,l} X_{i,k}  X_{j,m}' = - \frac{q^2\bw_{i,k} - \bw_{i,l}}{q^2\bw_{i,l} - \bw_{i,k}} X_{i,k} X_{i,l}  X_{j,m}', 
\]
it follows that $\Sym_{k,l} [X_{i,k},[X_{i,l}, X'_{j,m}]_{q\mi}]_q = 0$, and so $R_{k,l,m} + R_{l,k,m} = 0$. Next, 
\eqref{eq: Serre aux3}--\eqref{eq: Serre aux4} imply that 
\begin{align}
\gamma''_{l,k,m} &= \frac{q^{-2}(q^4C\mi\bw_{i,l}\mi - \bw_{i,k})}{\bw_{i,k} - C\mi\bw_{i,l}\mi} \gamma'_{k,l,m}, \\ 
X'_{i,l} X_{i,k}  X_{j,m} &= - \frac{q^2(\bw_{i,k}-C\mi \bw_{i,l}\mi)}{q^4C\mi\bw_{i,l}\mi-\bw_{i,k}} X_{i,k} X'_{i,l}  X_{j,m}. 
\end{align} 
Hence $[X_{i,k},[X'_{i,l}, X_{j,m}]_{q\mi}]_q + [X'_{i,l},[X_{i,k}, X_{j,m}]_{q\mi}]_q= 0$, and so $S_{k,l,m} + T_{l,k,m} = 0$. In summary, 
\eq \label{eq: MN vanish2} 
R_{k,l,m} + R_{l,k,m} = S_{k,l,m} + T_{l,k,m} = 0. 
\eneq

Define $P^*_{k,m}$, $R_{k,l,m}^*$, etc., by replacing $X \leftrightarrow X'$ and $\delta(a/b) \leftrightarrow \delta(aq^{-2} C b)$ (for $a \in \{w_1, w_2, z\}$ and $b \in \{\bw_{i,k}, \bw_{i,l}, \bw_{j,m}\}$) in the defining expression for $P_{k,m}$, $R_{k,l,m}$, etc., respectively. Then an analogous argument to the one above shows that 
\eq \label{eq: MN vanish3} 
M_{k,m}^* = N_{k,l,m}^* + N_{l,k,m}^* = P^*_{k,m} = R^*_{k,l,m} + R^*_{l,k,m} = S^*_{k,l,m} + T^*_{l,k,m} = 0. 
\eneq

\subsubsection{The vanishing coefficients: $\theta_i = 0, \theta_j = 1$} 
Since the argument is very similar to the one above, we will restrict ourselves to pointing out the main differences. 

\Lem \label{lem:note-theta-j-one}
Let $1\leq k\neq l\leq m_i$.
Then
\begin{align}
[X_{i,k},[X_{i,k},X''_j]_{q^{-1}}]_q&=0,
\label{eq:note-jpp-1}\\[1pt]
[X_{i,k},[X_{i,l},X''_j]_{q^{-1}}]_q
&=\widehat\gamma_{k,l}X_{i,k}X_{i,l}X''_j,
\label{eq:note-jpp-2}\\[1pt]
[X_{i,k},[X'_{i,l},X''_j]_{q^{-1}}]_q
&=\widehat\gamma'_{k,l}X_{i,k}X'_{i,l}X''_j,
\label{eq:note-jpp-3}\\[1pt]
[X'_{i,k},[X_{i,l},X''_j]_{q^{-1}}]_q
&=\widehat\gamma''_{k,l}X'_{i,k}X_{i,l}X''_j, 
\label{eq:note-jpp-4}
\end{align}
where 
\begin{align}
\widehat\gamma_{k,l}
&=
-\frac{
q^2(q-q^{-1})^2\bw_{i,k}\bw_{i,l}
\bigl(
\bw_{i,k}+\bw_{i,l}
-(q+q^{-1})C^{-\frac12}
\bigr)
}{
\bigl(q\bw_{i,k}-C^{-\frac12}\bigr)
\bigl(q\bw_{i,l}-C^{-\frac12}\bigr)
\bigl(\bw_{i,k}-q^2\bw_{i,l}\bigr)
},
\label{eq:note-hat-gamma}\\
\widehat\gamma'_{k,l}
&=
\frac{
q^4C^{-\frac12}\bw_{i,k}(q-q^{-1})^2
\bigl(
\bw_{i,k}\bw_{i,l}
+q^2C^{-1}
-(q+q^{-1})C^{-\frac12}\bw_{i,l}
\bigr)
}{
\bigl(q\bw_{i,k}-C^{-\frac12}\bigr)
\bigl(\bw_{i,l}-q^3C^{-\frac12}\bigr)
\bigl(\bw_{i,k}\bw_{i,l}-q^4C^{-1}\bigr)
},
\label{eq:note-hat-gamma-prime}\\
\widehat\gamma''_{k,l}
&=
-\frac{
q^2C^{-\frac12}\bw_{i,l}(q-q^{-1})^2
\bigl(
\bw_{i,k}\bw_{i,l}
+q^2C^{-1}
-(q+q^{-1})C^{-\frac12}\bw_{i,k}
\bigr)
}{
\bigl(q\bw_{i,l}-C^{-\frac12}\bigr)
\bigl(\bw_{i,k}\bw_{i,l}-C^{-1}\bigr)
\bigl(\bw_{i,k}-q^3C^{-\frac12}\bigr)
}.
\label{eq:note-hat-gamma-bis}
\end{align}
Moreover,
\begin{align}
\widehat\gamma_{k,l}
+q^2\frac{\bw_{i,k}-q^{-2}\bw_{i,l}}
 {\bw_{i,k}-q^2\bw_{i,l}} \widehat\gamma_{l,k}=0, \qquad 
\widehat\gamma'_{k,l}
+q^2\frac{\bw_{i,k}-C^{-1}\bw_{i,l}^{-1}}
 {\bw_{i,k}-q^4C^{-1}\bw_{i,l}^{-1}}\widehat\gamma''_{l,k}=0.
\end{align} 
\enlem 

\Proof
Notice that 
\[
\widehat\gamma_{k,l}
=\left.\gamma_{k,l,m}\right|_{\bw_{j,m}=q^2C\mhlf},
\qquad
\widehat\gamma'_{k,l}
=\left.\gamma'_{k,l,m}\right|_{\bw_{j,m}=C\mhlf},
\qquad
\widehat\gamma''_{k,l}
=\left.\gamma''_{k,l,m}\right|_{\bw_{j,m}=C\mhlf}.
\] 
The proof is formally the same as that of Lemma \ref{lem: XXX commutators}, with appropriate substitutions throughout. 
\enproof

For $k \neq l$, let  
\begin{align*} 
\widehat{P}_{k} :=& \ \Sym_{w_1, w_2} [\delta(w_1/ \bw_{i,k}) X_{i,k},[\delta(w_2/ \bw_{i,k}) X_{i,k}, \delta(zC\hlf)  X''_{j}]_{q\mi}]_q  \\ 
=& \  [X_{i,k},[X_{i,k}, X''_{j}]_{q\mi}]_q\Sym_{w_1, w_2} \delta(w_1/\bw_{i,k}) \delta(w_2/q^2\bw_{i,k})  \delta(zC\hlf)  , \\ 
\widehat{R}_{k,l} :=& \ \Sym_{w_1, w_2} [\delta(w_1/ \bw_{i,k}) X_{i,k},[\delta(w_2/ \bw_{i,l}) X_{i,l}, \delta(zC\hlf)  X''_{j}]_{q\mi}]_q \\ 
=& \ [X_{i,k},[X_{i,l}, X''_{j}]_{q\mi}]_q\Sym_{w_1, w_2}  \delta(w_1/\bw_{i,k}) \delta(w_2/\bw_{i,l})  \delta(zC\hlf) , \\ 
\widehat{S}_{k,l} :=& \ \Sym_{w_1, w_2} [\delta(w_1/ \bw_{i,k}) X_{i,k},[\delta(w_2 q^{-2} C \bw_{i,l}) X'_{i,l}, \delta(zC\hlf)  X''_{j}]_{q\mi}]_q \\ 
=& \  [X_{i,k},[X'_{i,l}, X''_{j}]_{q\mi}]_q\Sym_{w_1, w_2}  \delta(w_1/\bw_{i,k}) \delta(w_2q^{-2}C\bw_{i,l})  \delta(zC\hlf). 
\end{align*} 
Lemma \ref{lem:note-theta-j-one} (and its $*$-version) implies that 
\[
\widehat{P}_{k} = \widehat{P}_{k}^* = \widehat{R}_{k,l} + \widehat{R}_{l,k} = \widehat{R}^*_{k,l} + \widehat{R}^*_{l,k} = S_{k,l} + S^*_{l,k} = 0.  
\]

\subsubsection{The vanishing coefficients: $\theta_i = 1, \theta_j = 0$} 
We again restrict ourselves to pointing out the main differences. 

\Lem \label{lem:note-theta-i-one}
Let $1\leq k\leq m_i$ and $1 \leq m \leq m_j$. Then 
\begin{align}
[X''_i,[X_{i,k},X_{j,m}]_{q\mi}]_q
&=\breve{\gamma}_{k,m}X_{i,k}X''_iX_{j,m},\\
[X_{i,k},[X''_i,X_{j,m}]_{q\mi}]_q
&=-\breve{\gamma}_{k,m}X_{i,k}X''_iX_{j,m}, \\
[X''_i,[X_{i,k},X'_{j,m}]_{q\mi}]_q
&=\breve{\gamma}'_{k,m}X_{i,k}X''_iX'_{j,m},\\
[X_{i,k},[X''_i,X'_{j,m}]_{q\mi}]_q
&=-\breve{\gamma}'_{k,m}X_{i,k}X''_iX'_{j,m},
\end{align}
where
\begin{align}
\breve{\gamma}_{k,m}
&=-q^2C\mhlf
\frac{\bw_{i,k}(q-q\mi)^2\bigl( \bw_{i,k}+C\mhlf-(q+q\mi)\bw_{j,m}\bigr)}
{(\bw_{j,m}-qC\mhlf)(\bw_{i,k}-q^2C\mhlf)
 (q\bw_{i,k}-\bw_{j,m})}, \\
\breve{\gamma}'_{k,m}
&=(q-q\mi)^2
\frac{\bw_{i,k}\bw_{j,m}
\bigl(\bw_{i,k}\bw_{j,m}+C\mhlf\bw_{j,m}-(q^3+q)C\mi\bigr)}
{(\bw_{j,m}-qC\mhlf)(\bw_{i,k}-q^2C\mhlf)
 (\bw_{i,k}\bw_{j,m}-qC\mi)}. 
\end{align}
\enlem

\Proof
By Lemma \ref{lem: XiXj},
\[
X''_iX_{i,k}=rX_{i,k}X''_i,\qquad
X''_iX_{j,m}=aX_{j,m}X''_i,\qquad
X''_iX'_{j,m}=bX'_{j,m}X''_i, 
\]
for 
\[
r=\frac{q^2\bw_{i,k}-C^{-\frac12}}
{\bw_{i,k}-q^2C^{-\frac12}},\qquad
a=\frac{\bw_{j,m}-qC^{-\frac12}}
{q\bw_{j,m}-C^{-\frac12}},\qquad
b=q\frac{\bw_{j,m}-qC^{-\frac12}}
{\bw_{j,m}-q^3C^{-\frac12}}.
\]
and 
\[
X_{j,m}X_{i,k}=\beta X_{i,k}X_{j,m},\qquad
X'_{j,m}X_{i,k}=\beta'X_{i,k}X'_{j,m}, 
\]
for 
\[
\beta=q^{-1}\frac{\bw_{i,k}-q\bw_{j,m}}
{\bw_{i,k}-q^{-1}\bw_{j,m}},\qquad
\beta'=q^{-1}\frac{\bw_{i,k}-q^3C^{-1}\bw_{j,m}^{-1}}
{\bw_{i,k}-qC^{-1}\bw_{j,m}^{-1}}.
\]
A direct simplification gives
\[
r(1-q\mi\beta')+b\mi(\beta'-q)=\breve\gamma'_{k,m}, \qquad
1-q\mi b\mi-qr\beta'+r\beta'b\mi=-\breve\gamma'_{k,m}.
\]
This proves the last two equalities. For the first two equalities, the same calculation with $(b,\beta')$ replaced by
$(a,\beta)$ gives
\begin{align*}
[X''_i,[X_{i,k},X_{j,m}]_{q\mi}]_q
&=\Bigl(r(1-q\mi\beta)+a\mi(\beta-q)\Bigr)
 X_{i,k}X''_iX_{j,m},\\
[X_{i,k},[X''_i,X_{j,m}]_{q\mi}]_q
&=\Bigl(1-q\mi a\mi-qr\beta+r\beta a\mi\Bigr)
 X_{i,k}X''_iX_{j,m}.
\end{align*}
The two scalar coefficients simplify respectively to $\breve\gamma_{k,m}$ and
$-\breve\gamma_{k,m}$. 
\enproof

Let
\begin{align*} 
\breve{P}_{k,m} :=& \ \Sym_{w_1, w_2} [\delta(w_1C\hlf) X''_{i},[\delta(w_2/ \bw_{i,k}) X_{i,k}, \delta(zq^{-2}C\bw_{j,m})  X'_{j,m}]_{q\mi}]_q \\ 
=& \ [X''_{i},[X_{i,k}, X'_{j,m}]_{q\mi}]_q\Sym_{w_1, w_2}  \delta(w_1C\hlf) \delta(w_2/\bw_{i,k})  \delta(zq^{-2}C\bw_{j,m}) , \\ 
\breve{R}_{k,m} :=& \ \Sym_{w_1, w_2} [\delta(w_1/ \bw_{i,k}) X_{i,k},[\delta(w_2C\hlf) X''_{i}, \delta(zq^{-2}C\bw_{j,m})  X'_{j,m}]_{q\mi}]_q \\ 
=& \ [X_{i,k},[X''_{i}, X'_{j,m}]_{q\mi}]_q\Sym_{w_1, w_2}  \delta(w_1/\bw_{i,k}) \delta(w_2C\hlf)  \delta(zq^{-2}C\bw_{j,m}) , \\ 
\breve{S}_{k,m} :=& \ \Sym_{w_1, w_2} [\delta(w_1/ \bw_{i,k}) X_{i,k},[\delta(w_2 C\hlf) X''_{i}, \delta(z/\bw_{j,m})  X_{j,m}]_{q\mi}]_q \\ 
=& \  [X_{i,k},[X''_{i}, X_{j,m}]_{q\mi}]_q\Sym_{w_1, w_2}  \delta(w_1/\bw_{i,k}) \delta(w_2C\hlf)  \delta(z/\bw_{j,m}), \\
\breve{T}_{k,m} :=& \ \Sym_{w_1, w_2} [\delta(w_1C\hlf) X''_{i},[\delta(w_2/\bw_{i,k}) X_{i,k}, \delta(z/\bw_{j,m})  X_{j,m}]_{q\mi}]_q \\ 
=& \  [X''_{i},[X_{i,k}, X_{j,m}]_{q\mi}]_q\Sym_{w_1, w_2}  \delta(w_1C\hlf) \delta(w_2/\bw_{i,k})  \delta(z/\bw_{j,m}). 
\end{align*} 
Lemma \ref{lem:note-theta-i-one} (and its $*$-version) implies that 
\[
\breve{P}_{k,m} + \breve{R}_{k,m} = \breve{S}_{k,m} + \breve{T}_{k,m} = \breve{P}^*_{k,m} + \breve{R}^*_{k,m} = \breve{S}^*_{k,m} + \breve{T}^*_{k,m} = 0.  
\]

\subsubsection{All cases together} 

Finally, observe that the LHS of \eqref{eq: drel6} can be expressed as 
\begin{align}
\Sym _{w_1,w_2} &\big( \bA_i(w_1)\bA_i(w_2)\bA_{j}(z) - [2] \bA_i(w_1)\bA_{j}(z)\bA_i(w_2)+\bA_{j}(z)\bA_i(w_1)\bA_i(w_2) \big) = \\
&= \Sym_{w_1,w_2}[ A_i(w_1), [A_i(w_2), A_j(z)]_{q\mi}]_q.
\end{align} 
Abbreviate $\zeta_1 = w_1^{-\mu_i} \kappa_i \rho'_i$, $\zeta_2 = w_2^{-\mu_i} \kappa_i \rho'_i$ and $\zeta_3 = z^{-\mu_j} \kappa_j \rho'_j$.  
Putting together the results of the preceding subsections, we get 
\begin{align}
\TV &= \sum_{m = 1}^{m_j} \Big( \sum_{k=1}^{m_i} ( \zeta_1\zeta_2\zeta_3 M_{k,m} +  M_{k,m}^* )  + \sum_{1 \leq k \neq l \leq m_i} (  \zeta_1\zeta_2\zeta_3  N_{k,l,m} +  N_{k,l,m}^*)  \\
& \quad + \sum_{k=1}^{m_i} ( \zeta_1\zeta_2 P_{k,m} + \zeta_3 P_{k,m}^* ) + \sum_{1 \leq k \neq l \leq m_i} ( \zeta_1\zeta_2 R_{k,l,m} + \zeta_3 R_{k,l,m}^*)  \\
& \quad + \sum_{1 \leq k \neq l \leq m_i} ( \zeta_1 \zeta_3 S_{k,l,m} +  \zeta_2 S_{k,l,m}^*) + \sum_{1 \leq k \neq l \leq m_i} (\zeta_2\zeta_3 T_{k,l,m} +  \zeta_1 T_{k,l,m}^* ) \Big) \\ 
& \quad + \theta_j \tau_j \Big(\sum_{k=1}^{m_i}  (  \zeta_1\zeta_2 \widehat{P}_k + \widehat{P}^*_k) + \sum_{1 \leq k \neq l \leq m_i}  \big(  \zeta_1\zeta_2 \widehat{R}_{k,l} + \widehat{R}^*_{k,l} +  \zeta_1\widehat{S}_{k,l} + \zeta_2\widehat{S}^*_{k,l} \big) \Big)   \\ 
& \quad + \theta_i \tau_i \sum_{m=1}^{m_j} \sum_{k=1}^{m_i} \Big( 
 \zeta_2 \breve{P}_{k,m}+ \zeta_1\breve{R}_{k,m} + \zeta_3(\breve{P}^*_{k,m} + \breve{R}^*_{k,m}) \\ 
& \quad \quad +  \zeta_1\zeta_3\breve{S}_{k,m}+ \zeta_2\zeta_3\breve{T}_{k,m} + \breve{S}^*_{k,m} + \breve{T}^*_{k,m}
\Big) \\
&= 0.
\end{align}

\subsubsection{The remainder} 

Let us denote the remainder which does not vanish on the LHS as $\mathbf{Q}$. We do not make any assumptions about $\theta_i, \theta_j$, except the standard constraint $\theta_i \theta_j = 0$. Then, by definition, 
\begin{align}
\mathbf{Q} &= \kappa_i \rho'_i \sum_{k=1}^{m_i} \Sym_{w_1,w_2} w_1^{-\mu_i} \Big( 
\mathbf{F} ( \delta(w_1/\bw_{i,k}) X_{i,k}, \delta(w_2q^{-2}C\bw_{i,k}) X'_{i,k}, \bA_j(z) ) \\ 
&\quad + \mathbf{F} (\delta(w_2q^{-2}C\bw_{i,k}) X'_{i,k}, \delta(w_1/\bw_{i,k}) X_{i,k}, \bA_j(z) ) \Big) \\
&\quad + \theta_i 2\tau_i^2 \delta(w_1C\hlf)\delta(w_2 C\hlf)  \mathbf{F}(X''_{i}, X''_{i}, \bA_j(z)),  
\end{align}
where 
\[
\mathbf{F}(a,b,c) = abc -[2]acb + cab. 
\]
We will use the following lemma.

\Lem \label{lem: XXX'=W}
Let $\theta_j = 0$. Then the following identities hold: 
\begin{align*} 
\mathbf{F}(X_{i,k}, X'_{i,k}, Y)    &=  \frac{\rho\Omega \mathfrak{W}_{i,k}(\bw_{i,k})}{(q^{-1} \bw_{i,k} {-} q^{-2} \bw_{j,m}) ( q^{-1} \bw_{i,k} {-} C^{-1} \bw_{j,m}^{-1} ) } 
 Y, \\
\mathbf{F}(X'_{i,k}, X_{i,k}, Y)    &= \frac{-\rho q^2\Omega\mathfrak{W}_{i,k}(q^{-2}\bw_{i,k})}{(q^{-3} \bw_{i,k} {-} q^{-2} \bw_{j,m}) ( q^{-3} \bw_{i,k} {-} C^{-1} \bw_{j,m}^{-1} )} 
 Y, \\ 
\theta_i \mathbf{F}(X''_{i}, X''_{i}, Y)    &=  \theta_i \frac{ (-1)^{\mu_i}C^{\frac{\mu_i-1}{2}} q^{3(\mu_i - \lambda_i)}(q^2-1)^2\bw_{j,m} q\mi \mathfrak{W}_{i}(C\mhlf)}{(\bw_{j,m} - qC\mhlf)^2}Y,
\end{align*}
where 
\[
\Omega =  \frac{1}{\bw_{i,k} (q^{-2}\bw_{i,k} {-} C\mi \bw_{i,k}\mi)} 
\] 
and $Y \in \{X_{j,m}, X'_{j,m}\}_{m=1}^{m_j}$. Next, let $\theta_j = 1$. Then 
\begin{align}
\mathbf{F}(X_{i,k},X'_{i,k},X''_j)
&=\frac{\rho\Omega\mathfrak{W}_{i,k}(\bw_{i,k})}
{(q\mi\bw_{i,k}-q^{-2}C\mhlf)(q\mi\bw_{i,k}-C\mhlf)}X''_j,
\label{eq:note-F-Xpp-1}\\
\mathbf{F}(X'_{i,k},X_{i,k},X''_j)
&=\frac{-\rho q^2\Omega\mathfrak{W}_{i,k}(q^{-2}\bw_{i,k})}
{(q^{-3}\bw_{i,k}-q^{-2}C\mhlf)(q^{-3}\bw_{i,k}-C\mhlf)}X''_j.
\label{eq:note-F-Xpp-2}
\end{align}
\enlem 

\Proof
Let us verify the first identity, leaving the others to the reader. Indeed, we get
\begin{align*}
& X_{i,k}X_{i,k}'X_{j,m} - [2] X_{i,k}X_{j,m}X_{i,k}' + X_{j,m}X_{i,k}X_{i,k}' = \\
& =  X_{i,k}X_{i,k}'X_{j,m} - [2] X_{i,k}  q \frac{ \bw_{j,m} - q C^{-1} \bw_{i,k}^{-1} }{ \bw_{j,m} - q^3 C^{-1} \bw_{i,k}^{-1} } X_{i,k}' X_{j,m}\\
& \quad + q^{-1} \frac{\bw_{i,k} - q \bw_{j,m}}{\bw_{i,k} - q^{-1} \bw_{j,m}} X_{i,k}  q \frac{ \bw_{j,m} - q C^{-1} \bw_{i,k}^{-1} }{ \bw_{j,m} - q^3 C^{-1} \bw_{i,k}^{-1} }   X_{i,k}'X_{j,m} \\ 
& =  (1 -[2]q \frac{ \bw_{j,m} - q^{-1} C^{-1} \bw_{i,k}^{-1} }{ \bw_{j,m} - q C^{-1} \bw_{i,k}^{-1} } + \frac{ \bw_{j,m} - q^{-1} C^{-1} \bw_{i,k}^{-1} }{ \bw_{j,m} - q C^{-1} \bw_{i,k}^{-1} }\frac{\bw_{i,k} - q \bw_{j,m}}{\bw_{i,k} - q^{-1} \bw_{j,m}}  ) \times \\
& \quad \times X_{i,k}X_{i,k}'X_{j,m} \\ 
&= \frac{-q^{-1}(q-q\mi)\bw_{i,k}(\bw_{i,k} - q^{-2}C\mi \bw_{i,k}\mi)}{( q^{-1} \bw_{i,k} {-} C^{-1} \bw_{j,m}^{-1} )(q\mi\bw_{i,k} - q^{-2} \bw_{j,m})} X_{i,k}X_{i,k}'X_{j,m} \\ 
&= \frac{\rho\Omega \mathfrak{W}_{i,k}(\bw_{i,k})}{(q^{-1} \bw_{i,k} {-} q^{-2} \bw_{j,m}) ( q^{-1} \bw_{i,k} {-} C^{-1} \bw_{j,m}^{-1} ) } 
 X_{j,m}. 
\end{align*} 
The first equality above follows from Lemma \ref{lem: XiXj}, and the fourth from~Lemma~\ref{lem: XX' equal}. 
\enproof

Let us redefine the symbols from \eqref{eq: deltas1}--\eqref{eq: deltas2} in the following fashion: 
\begin{align*}
\delta_1 &= \delta(w_1/\bw_{i,k}) \delta(w_2 C \bw_{i,k}), \quad& 
\delta_2 &= \delta(q^2w_2/\bw_{i,k}) \delta(w_1 q^{-2} C \bw_{i,k}),  \\ 
\delta_3 &= \delta(w_2/\bw_{i,k}) \delta(w_1  C \bw_{i,k}), \quad&  
\delta_4 &=  \delta(q^2w_1/\bw_{i,k}) \delta(w_2 q^{-2} C \bw_{i,k}). 
\end{align*} 
\[
\delta_5 = \delta(w_1 C\hlf) \delta(w_2 C\hlf). 
\]
It follows from Lemma \ref{lem: XXX'=W} that 
\begin{align} 
\mathbf{Q} &=  \kappa_i \rho'_i \sum_{k=1}^{m_i} \Sym_{w_1,w_2} w_1^{-\mu_i} \Big( 
\delta(w_1/\bw_{i,k})\delta(w_2C\bw_{i,k}) \mathbf{F}(X_{i,k}, X'_{i,k}, \As{j}{z}) \\ 
&\quad + \delta(q^2w_1/\bw_{i,k})\delta(w_2Cq^{-2}\bw_{i,k}) \mathbf{F}(X'_{i,k}, X_{i,k}, \As{j}{z}) \Big) \\ 
&\quad + \theta_i 2\tau_i^2 \delta_5  \mathbf{F}(X''_{i}, X''_{i}, \bA_j(z)) \\
 \label{eq: Q exp}
&= \kappa_i \rho'_i \rho \sum_{k=1}^{m_i} \Big(  \frac{\bw_{i,k}^{-1-\mu_i}(\delta_1 +  \delta_3)\mathfrak{W}_{i,k}(\bw_{i,k}) }{(q^{-2}\bw_{i,k} {-} C\mi \bw_{i,k}\mi)(q^{-1} \bw_{i,k} {-} q^{-2} z) ( q^{-1} \bw_{i,k} {-} C^{-1} z^{-1} ) } 
\As{j}{z} \\ 
&\quad -  \frac{q^{2+2\mu_i} \bw_{i,k}^{-1-\mu_i}(\delta_2 + \delta_4) \mathfrak{W}_{i,k}(q^{-2}\bw_{i,k}) }{(q^{-2}\bw_{i,k} {-} C\mi \bw_{i,k}\mi) (q^{-3} \bw_{i,k} {-} q^{-2} z) ( q^{-3} \bw_{i,k} {-} C^{-1} z^{-1} )} 
 \As{j}{z} \Big) \\ 
 &\quad + \theta_i \delta_5 \frac{2q\mi C^{\frac{\mu_i-1}{2}}(1-q)^2 \kappa_i \rho'_i z \mathfrak{W}_i(C\mhlf)}{(z-qC\mhlf)^2}\As{j}{z}. 
\end{align}

\subsection{The RHS} 

If we write 
\[
\mathbf{G}(a) = [2]zw_2C [ a, \As{j}{z}]_{q^{-2}} + (1+Cw_2^2) [ \As{j}{z}, a]_{q^{-2}}, 
\]
then Lemma \ref{lem: DeltaThetaX} implies that the RHS of \eqref{eq: drel6} equals 
\begin{align} 
\RHS &= \frac{- \kappa_i \rho \bDel(w_1w_2) (1-Cw_2^2)}{(1-q^2Cw_2^2)(1-q^{-2}Cw_2^2)} 
\mathbf{G}(\aThs{i}{w_2} - \aThs{i}{C\mi w_2\mi}) \\
&= \frac{- w_2 C}{(1-q^2Cw_2^2)(1-q^{-2}Cw_2^2)} \mathbf{G}(\mathfrak{X}_i) \\ \label{eq: RHS Serre pr}
&=  \frac{w_2 C \kappa_i \rho'_i \rho^2}{(1-q^2Cw_2^2)(1-q^{-2}Cw_2^2)} \sum_{k=1}^{m_i} 
 \bw_{i,k}^{-2-\mu_i}(q^{-2}\bw_{i,k} {-} C\mi \bw_{i,k}\mi)\mi \times \\
&\quad \times \Big( q^3(\delta_3 + \delta_1) \mathbf{G}(\mathfrak{W}_{i,k}(\bw_{i,k})) - q^{7+2\mu_i}(\delta_4 + \delta_2) \mathbf{G}( \mathfrak{W}_{i,k}(q^{-2}\bw_{i,k})) \Big) \\ 
&\quad - \theta_i \frac{ 2\kappa_i \rho_i' (1-q) C^{\frac{\mu_i-1}{2}}  C\hlf  }{(1+q)(q-q\mi)^2}\delta_5\mathbf{G}(\mathfrak{W}_i(C\mhlf)). 
\end{align} 
A short calculation using Lemma \ref{lem: XjTi} shows that 
\begin{align*}
\delta_1\mathbf{G}(\mathfrak{W}_{i,k}(\bw_{i,k})) &= \frac{(q - q\mi)(1 - q^{-2}Cw_2^{2})(1-Cq^{2}w_2^{2})}{C^{2}q^3w_2^{2}(q^{-1} \bw_{i,k} {-} q^{-2} z) ( q^{-1} \bw_{i,k} {-} C^{-1} z^{-1} )} \delta_1 \mathfrak{W}_{i,k}(\bw_{i,k}) \As{j}{z}, \\
\delta_4\mathbf{G}(\mathfrak{W}_{i,k}(q^{-2}\bw_{i,k})) &= \frac{(q - q\mi)(1 - q^{-2}Cw_2^{2})(1-Cq^{2}w_2^{2})}{C^{2}q^3w_2^{2}(q^{-3} \bw_{i,k} {-} q^{-2} z) ( q^{-3} \bw_{i,k} {-} C^{-1} z^{-1} )} \delta_4 \mathfrak{W}_{i,k}(q^{-2}\bw_{i,k}) \As{j}{z}, \\ 
\delta_3\mathbf{G}(\mathfrak{W}_{i,k}(\bw_{i,k})) &= \frac{(q - q\mi)(1 - q^{-2}Cw_2^{2})(1-Cq^{2}w_2^{2})}{Cq^3(q^{-1} \bw_{i,k} {-} q^{-2} z) ( q^{-1} \bw_{i,k} {-} C^{-1} z^{-1} )} \delta_3 \mathfrak{W}_{i,k}(\bw_{i,k}) \As{j}{z}, \\ 
\delta_2\mathbf{G}(\mathfrak{W}_{i,k}(q^{-2}\bw_{i,k})) &= \frac{(q - q\mi)(1 - q^{-2}Cw_2^{2})(1-Cq^{2}w_2^{2})}{Cq^3(q^{-3} \bw_{i,k} {-} q^{-2} z) ( q^{-3} \bw_{i,k} {-} C^{-1} z^{-1} )} \delta_2 \mathfrak{W}_{i,k}(q^{-2}\bw_{i,k}) \As{j}{z}, \\
\delta_5 \mathbf{G}(\mathfrak{W}_i(C\mhlf)) &= \frac{C\mhlf(q-q\mi)^3z}{(z-qC\mhlf)^2} \delta_5 \mathfrak{W}_i(C\mhlf) \mathbf{A}_j(z). 
\end{align*}
Therefore, 
\begin{align}
\frac{w_2 C \kappa_i \rho'_i \rho^2  \bw_{i,k}^{-2-\mu_i} q^3 \delta_1}{(1-q^2Cw_2^2)(1-q^{-2}Cw_2^2)} &\mathbf{G}(\mathfrak{W}_{i,k}(\bw_{i,k})) = \\ \label{eq: G-A1} 
 &= \frac{\kappa_i \rho'_i \rho  \bw_{i,k}^{-1-\mu_i}\delta_1 \mathfrak{W}_{i,k}(\bw_{i,k}) }{(q^{-1} \bw_{i,k} {-} q^{-2} z) ( q^{-1} \bw_{i,k} {-} C^{-1} z^{-1} )} \As{j}{z}, \\ 
\frac{w_2 C \kappa_i \rho'_i \rho^2  \bw_{i,k}^{-2-\mu_i} q^{7+2\mu_i} \delta_4}{(1-q^2Cw_2^2)(1-q^{-2}Cw_2^2)} &\mathbf{G}(\mathfrak{W}_{i,k}(q^{-2}\bw_{i,k})) = \\ \label{eq: G-A2} \tag{\theequation} \stepcounter{equation}
 &= \frac{\kappa_i \rho'_i \rho  \bw_{i,k}^{-1-\mu_i}q^{2+2\mu_i}\delta_4 \mathfrak{W}_{i,k}(q^{-2}\bw_{i,k})}{(q^{-3} \bw_{i,k} {-} q^{-2} z) ( q^{-3} \bw_{i,k} {-} C^{-1} z^{-1} )}  \As{j}{z}, \\ 
\frac{w_2 C \kappa_i \rho'_i \rho^2  \bw_{i,k}^{-2-\mu_i} q^3 \delta_3}{(1-q^2Cw_2^2)(1-q^{-2}Cw_2^2)} &\mathbf{G}(\mathfrak{W}_{i,k}(\bw_{i,k})) = \\ \label{eq: G-A3} \tag{\theequation} \stepcounter{equation}
 &= \frac{\kappa_i \rho'_i \rho  \bw_{i,k}^{-1-\mu_i}\delta_3 \mathfrak{W}_{i,k}(\bw_{i,k}) }{(q^{-1} \bw_{i,k} {-} q^{-2} z) ( q^{-1} \bw_{i,k} {-} C^{-1} z^{-1} )} \As{j}{z}, \\ 
\frac{w_2 C \kappa_i \rho'_i \rho^2  \bw_{i,k}^{-2-\mu_i} q^{7+2\mu_i} \delta_2}{(1-q^2Cw_2^2)(1-q^{-2}Cw_2^2)} &\mathbf{G}(\mathfrak{W}_{i,k}(q^{-2}\bw_{i,k})) = \\ \label{eq: G-A4} \tag{\theequation} \stepcounter{equation}
 &= \frac{\kappa_i \rho'_i \rho  \bw_{i,k}^{-1-\mu_i}q^{2+2\mu_i}\delta_2 \mathfrak{W}_{i,k}(q^{-2}\bw_{i,k})}{(q^{-3} \bw_{i,k} {-} q^{-2} z) ( q^{-3} \bw_{i,k} {-} C^{-1} z^{-1} )}  \As{j}{z}, \\
  -\theta_i \frac{ 2\kappa_i \rho_i' (1-q) C^{\frac{\mu_i-1}{2}}  C\hlf  }{(1+q)(q-q\mi)^2}\delta_5 &\mathbf{G}(\mathfrak{W}_i(C\mhlf)) = \\  \label{eq: G-A5}
  &= \theta_i \frac{2q\mi C^{\frac{\mu_i-1}{2}}(1-q)^2 \kappa_i \rho'_i z}{(z-qC\mhlf)^2}  \delta_5 \mathfrak{W}_i(C\mhlf) \mathbf{A}_j(z). 
\end{align} 
Comparing \eqref{eq: Q exp} with \eqref{eq: RHS Serre pr} (using \eqref{eq: G-A1}--\eqref{eq: G-A5}), we conclude that $\mathbf{Q} = \RHS$.


\addtocontents{toc}{\SkipTocEntry}
\section*{Conflict of interest statement} 

On behalf of all authors, the corresponding author states that there is no conflict of interest. 

\addtocontents{toc}{\SkipTocEntry}
\section*{Data availability statement} 

The authors declare that all the data supporting the findings of this article are available within the paper. 

\providecommand{\bysame}{\leavevmode\hbox to3em{\hrulefill}\thinspace}
\providecommand{\MR}{\relax\ifhmode\unskip\space\fi MR }
\providecommand{\MRhref}[2]{%
  \href{http://www.ams.org/mathscinet-getitem?mr=#1}{#2}
}
\providecommand{\href}[2]{#2}


\begin{thebibliography}{KWWY14}

\bibitem[AV22]{AppelVlaar}
A. Appel, B. Vlaar, \emph{Universal {K}-matrices for quantum {K}ac-{M}oody
  algebras}, Represent. Theory \textbf{26} (2022), 764--824. \MR{4454332}

\bibitem[AV25a]{AppelVlaar4}
A. Appel, B. Vlaar, \emph{Boundary transfer matrices arising from quantum symmetric
  pairs}, Indag. Math. (N.S.) \textbf{36} (2025), no.~6, 1830--1862.
  \MR{4978865}

\bibitem[AV25b]{AppelVlaar3}
A. Appel, B. Vlaar, \emph{Tensor {K}-{M}atrices for {Q}uantum {S}ymmetric {P}airs}, Comm.
  Math. Phys. \textbf{406} (2025), no.~5, Paper No. 100. \MR{4887619}

\bibitem[AV25c]{AppelVlaar2}
A. Appel, B. Vlaar, \emph{{Trigonometric K-matrices for finite-dimensional representations
  of quantum affine algebras}}, J. Eur. Math. Soc. (2025).

\bibitem[BFN19]{BFN}
A. Braverman, M. Finkelberg, H. Nakajima, \emph{Coulomb branches of {$3d$}
  {$\mathcal{N}=4$} quiver gauge theories and slices in the affine
  {G}rassmannian}, Adv. Theor. Math. Phys. \textbf{23} (2019), no.~1, 75--166,
  With two appendices by Braverman, Finkelberg, Joel Kamnitzer, Ryosuke Kodera,
  Nakajima, Ben Webster and Alex Weekes. \MR{4020310}

\bibitem[BK20]{bas-kol-20}
P. Baseilhac, S. Kolb, \emph{Braid group action and root vectors for the
  {$q$}-{O}nsager algebra}, Transform. Groups \textbf{25} (2020), no.~2,
  363--389. \MR{4098883}

\bibitem[BPT25]{BPT}
R. Bartlett, T. Prze\'{z}dziecki, L. Tappeiner, \emph{GKLO representations of
  twisted Yangians in type $\mathsf{AI}$ and quantizations of symmetric
  quotients of the affine Grassmannian}, 2025, \arxiv{2510.12706}.

\bibitem[FPT22]{FPT}
R. Frassek, V. Pestun, A. Tsymbaliuk, \emph{Lax matrices from antidominantly
  shifted {Y}angians and quantum affine algebras: {A}-type}, Adv. Math.
  \textbf{401} (2022), Paper No. 108283, 73. \MR{4394682}

\bibitem[FT19]{FT19}
M. Finkelberg, A. Tsymbaliuk, \emph{Multiplicative slices, relativistic {T}oda
  and shifted quantum affine algebras}, Representations and nilpotent orbits of
  {L}ie algebraic systems, Progr. Math., vol. 330, Birkh\"auser/Springer, Cham,
  2019, pp.~133--304. \MR{3971731}

\bibitem[GK91]{GavKl}
A.~M. Gavrilik, A.~U. Klimyk, \emph{{$q$}-deformed orthogonal and
  pseudo-orthogonal algebras and their representations}, Lett. Math. Phys.
  \textbf{21} (1991), no.~3, 215--220. \MR{1102131}

\bibitem[GKLO05a]{GKLO}
A. Gerasimov, S. Kharchev, D. Lebedev, S. Oblezin, \emph{On a class of
  representations of the {Y}angian and moduli space of monopoles}, Comm. Math.
  Phys. \textbf{260} (2005), no.~3, 511--525. \MR{2182434}

\bibitem[GKLO05b]{MR2184015}
A. Gerasimov, S. Kharchev, D. Lebedev, S. Oblezin, \emph{On a class of
  representations of quantum groups}, Noncommutative geometry and
  representation theory in mathematical physics, Contemp. Math., vol. 391,
  Amer. Math. Soc., Providence, RI, 2005, pp.~101--110. \MR{2184015}

\bibitem[KWWY14]{KWWY}
J. Kamnitzer, B. Webster, A. Weekes, O. Yacobi, \emph{Yangians and
  quantizations of slices in the affine {G}rassmannian}, Algebra Number Theory
  \textbf{8} (2014), no.~4, 857--893. \MR{3248988}

\bibitem[LP25a]{LP}
J.-R. Li, T. Prze\'{z}dziecki, \emph{Boundary $q$-characters of evaluation modules
  for split quantum affine symmetric pairs}, 2025, \arxiv{2504.14042}. 
  
\bibitem[LP25b]{LP25}
J.-R. Li, T. Prze\'{z}dziecki, \emph{{C}ompatibility of {D}rinfeld
  presentations for split affine {K}ac--{M}oody quantum symmetric pairs}, Lett.
  Math. Phys. \textbf{115} (2025), 75. 
  
\bibitem[LP26a]{LP26}
J.-R. Li, T. Prze\'{z}dziecki, \emph{Compatibility of Drinfeld presentations and $q$-characters for affine Kac-Moody quantum symmetric pairs: quasi-split case}, 2026, \arxiv{2601.02165}. 

\bibitem[LP26b]{iHopf} 
M. Lu and X. Pan, \emph{Drinfeld presentation of affine $\imath$quantum groups via $\imath$Hopf algebras}, preprint, 2026.

\bibitem[LPWZ25]{LPWZ}
M. Lu, X. Pan, W. Wang, W. Zhang, \emph{Braid group action and quasi-split
  affine iquantum groups {III}}, 2025, \arxiv{2511.00882}.

\bibitem[LRW23]{lu-ruan-wang-23}
M. Lu, S. Ruan, W. Wang, \emph{{$\imath$}{H}all algebra of the projective line
  and {$q$}-{O}nsager algebra}, Trans. Amer. Math. Soc. \textbf{376} (2023),
  no.~2, 1475--1505. \MR{4531682}

\bibitem[LRZ25]{MR4967120}
M. Lu, S. Ruan, H. Zhang, \emph{Analogue of {F}eigin's map on the
  {$\imath$}quantum groups of split type}, Int. Math. Res. Not. IMRN (2025),
  no.~19, Paper No. rnaf297, 28. \MR{4967120}

\bibitem[LW21]{lu-wang-21}
M. Lu, W. Wang, \emph{A {D}rinfeld type presentation of affine
  {$\imath$}quantum groups {I}: {S}plit {ADE} type}, Adv. Math. \textbf{393}
  (2021), Paper No. 108111, 46. \MR{4340233}

\bibitem[LWW25a]{LWW1}
K. Lu, W. Wang, A. Weekes, \emph{Shifted twisted {Y}angians and affine
  {G}rassmannian islices}, 2025, \arxiv{2510.10652}.

\bibitem[LWW25b]{LWW2}
K. Lu, W. Wang, A. Weekes, \emph{Shifted twisted {Y}angians of quasi-split {ADE} types}, 2025,
  \arxiv{2512.19998}. 
  
  \bibitem[LWW25c]{LWW3}
K. Lu, W. Wang, A. Weekes, \emph{Shifted affine $\imath$quantum groups of quasi-split {ADE} types}, 2026,
  \arxiv{2603.28446}.

\bibitem[LWZ23]{LWZ-I}
M. Lu, W. Wang, W. Zhang, \emph{Braid group action and quasi-split affine
  {$\imath$}quantum groups {I}}, Represent. Theory \textbf{27} (2023),
  1000--1040. \MR{4659754}

\bibitem[LWZ24]{LWZ-II}
M. Lu, W. Wang, W. Zhang, \emph{Braid {G}roup {A}ction and {Q}uasi-{S}plit {A}ffine {$\imath$}{Q}uantum {G}roups {II}: {H}igher {R}ank}, Comm. Math. Phys. \textbf{405}
  (2024), no.~6, Paper No. 142. \MR{4751724}

\bibitem[Nak25]{NakajimaQSP}
H. Nakajima, \emph{Instantons on ALE spaces for classical groups, involutions
  on quiver varieties, and quantum symmetric pairs}, 2025, \arxiv{2510.13007}. 

\bibitem[Prz25]{Przez-23}
T. Prze\'zdziecki, \emph{Drinfeld rational fractions for affine {K}ac-{M}oody
  quantum symmetric pairs}, Selecta Math. (N.S.) \textbf{31} (2025), no.~5,
  Paper No. 87, 59. \MR{4961836}

\bibitem[SSX25]{SSX}
Y. Shen, C. Su, R. Xiong, \emph{Quivers with involutions and shifted twisted
  {Y}angians via {C}oulomb branches}, 2025, \arxiv{2510.12118}.
  
\bibitem[SSX26]{Sunew}
Y. Shen, C. Su, R. Xiong, \emph{Quantized Coulomb Branches of Separated Cotangent Type and Orthosymplectic Quivers}, 
2026, \arxiv{2608.16091}. 

\bibitem[Zha22]{ZhangDr}
W. Zhang, \emph{A {D}rinfeld-type presentation of affine {$\imath$}quantum
  groups {II}: split {BCFG} type}, Lett. Math. Phys. \textbf{112} (2022),
  no.~5, Paper No. 89, 33. \MR{4482099}

\bibitem[Wan25]{Zichang}
Z. Wang, \emph{Quivers with involutions and shifted twisted Yangians via
  Coulomb branches II}, 2025, \arxiv{2601.00039}.

\end{thebibliography}
\end{document}